\magnification 1200
\input amssym.def
\input amssym.tex
\parindent = 40 pt
\parskip = 12 pt
\font \heading = cmbx10 at 12 true pt
 at 22 true pt
\font \medheading =cmbx7 at 16 true pt
 at 7 true pt
\def \R{{\bf R}}

\centerline{\medheading Uniform bounds for Fourier Transforms of Surface}
\centerline{\medheading Measures in R$^3$ with Nonsmooth Density}
\rm
\line{}
\line{}
\centerline{\heading Michael Greenblatt}
\line{}
\centerline{September 10, 2014}
\baselineskip = 12 pt
\font \heading = cmbx10 at 14 true pt
\line{}
\line{}
\line{}
\noindent{\heading 1. Background and Theorem Statements}.

\vfootnote{}{This research was supported in part by NSF grant  DMS-1001070} In this paper we consider Fourier transforms of measures
of the form $Q(s)d\mu$, where $d\mu$ denotes the surface measure corresponding to a bounded subset of a real-analytic hypersurface in $\R^3$ and $Q(s)$ denotes a function on the surface which may have singularities. To be precise, after a partition of unity and a translation and rotation of coordinates we consider oscillatory integrals of the form
$$T(\lambda_1,\lambda_2,\lambda_3) = \int_{\R^3}e^{i\lambda_1 S(x,y) + i\lambda_2 x + i\lambda_3 y}\,g(S(x,y))K(x,y)\,dx\,dy \eqno (1.1)$$
Here $S(x,y)$ is a nonconstant real-analytic phase function on a neighborhood of the origin such that $S(0,0) = 0$ and $\nabla S(0,0) = 0$. When viewed in terms of the hypersurface lying in $\R^3$, the density in $(1.1)$ is of the form $K(x,y)g(z)$. 

The functions $g(z)$ and $K(x,y)$ satisfy the following conditions.
The function $g(z)$ is assumed to be real-valued and $C^1$ on $\R - \{0\}$ such that for some real $\alpha$ and some $A > 0$ one has
$$|g(z)| \leq A|z|^{\alpha}\,\,\,\,\,\,\,\,\,\,\,\,\,\,\,\,\,\,\,\,\,\,\,\,\,\,\,\,|g'(z)| \leq A|z|^{\alpha - 1} \eqno (1.2a)$$
The function $K(x,y)$ is assumed to be a $C^1$ real-valued compactly supported function on $\R^2 - \{(0,0)\}$ such that 
for some  real $\beta$ we have
$$| K(x,y)| \leq  A(x^2 + y^2)^{\beta  \over 2}\,\,\,\,\,\,\,\,\,\,\,\,\,\,\,\,\,\,\,\,\,\,\,\,\,\,\,\,|\nabla  K(x,y)| \leq  A(x^2 + y^2)^{\beta   - 1\over 2}\eqno (1.2b)$$
Both $\alpha$ and $\beta$ can be negative, but to ensure that $(1.1)$ is well-defined we require that $(x^2 + y^2)^{\beta \over 2}|S(x,y)|^{\alpha}$ is integrable over some neighborhood of the origin.

In this paper, we will prove uniform estimates on $T(\lambda_1,\lambda_2,\lambda_3)$ that generalize the sharp uniform estimates that are known to hold in the situation where $K(x,y)g(z)$ is smooth function $\phi(x,y)$ (using $\alpha = \beta = 0$). The latter results can be described as follows.
By resolution of singularities (see Ch. 7 of [AGV] for details),  there are $\delta > 0$ and an integer $d = 0$ or $1$ such that if the support of $\phi(x,y)$ is sufficiently small, then for some $C_{\phi}$ as $\lambda_1 \rightarrow \infty$ one has
$$T(\lambda_1,0,0) = C_{\phi}\lambda_1^{-\delta} (\ln \lambda_1 )^d + o(|\lambda_1^{-\delta} (\ln \lambda_1 )^d|) \eqno (1.3)$$
Here $C_{\phi}$ will be nonzero if $\phi(x,y)$ is nonnegative with $\phi(0,0) > 0$. It is a consequence of the stability theorems of [K1][K2]  
that if the support of $\phi(x,y)$ is sufficiently small, then for some $C_{\phi}'$ one has uniform estimates 
$$|T(\lambda_1,\lambda_2,\lambda_3)| \leq C_{\phi}'|\lambda_1^{-\delta} (\ln \lambda_1 )^d| \eqno (1.4)$$
Analogous results for smooth $S(x,y)$ are shown in [D][IKeM2][IM2]. In this paper we will prove estimates of the form $(1.4)$ for
the possibly singular densities here (for real-analytic $S(x,y)$), with appropriately defined $\delta$ and $d$. These estimates will be uniform in $\lambda_2$ and $\lambda_3$, and will also be uniform over all densities satisfying $(1.2a)-(1.2b)$. They will imply the above results for smooth $\phi(x,y)$ whenever $S(x,y)$ has a zero of order greater than 5 at the origin, and for some of the situations where
it has a zero of order between 3 and 5 at the origin. The analysis of this paper is based on an explicit resolution of singularities algorithm
(Theorem 2.1) as well as versions of the Van der Corput lemma. We do not make use of the adapted coordinate systems that have
often been used in this subject.

It is worth pointing out that in the case of the oscillatory integral operators with phase functions of two variables that are considered in [PS1], extensions to weighted oscillatory integral operators are proven in [PS2]. Although there are a number of differences when dealing with operators, in [PS1]-[PS2] one also uses a type of resolution of singularities to divide into wedges on which one can do 
an appropriate analysis, and these papers influenced the development of the resolution of singularities algorithm used in this paper.

The question of proving optimal estimates for
two-dimensional oscillatory integrals given a fixed density function has been analyzed in [PrY]. There are also the recent papers [CKaN] [KaN], which also deal with proving optimal estimates for oscillatory integrals given a fixed smooth density function, for
classes of phases in any dimension. In addition, damped oscillatory integrals related to those of this paper often appear in the study of maximal averages over surfaces. We refer to [SoS] [IoSa1] [IoSa2] [IKeM1] for more on this. We also mention work by Lichtin [Lic] on related topics.

In the case of convex hypersurfaces in any dimension, for specific classes of densities connected to the surface, Fourier transform decay
estimates have been proven in several papers including [Lit][CoMa].

We turn to defining the $\delta$ and $d$ that we will use in our theorems. Let $E_r$ denote the disk $\{(x,y): x^2 + y^2 < r^2 \}$. In Lemma 2.2 we will show that there is a $\delta > 0$ and an integer $d = 0$ or $1$ such that if $r$ is sufficiently small 
there are constants $C_r$ and $C_r'$ such that for sufficiently small $\epsilon$ one has
$$ C_r \epsilon^{\delta}|\ln \epsilon|^d  \leq \int_{\{(x,y) \in E_r :  |S(x,y)| < \epsilon\}}|S(x,y)|^{ \alpha} (x^2 + y^2)^{\beta \over 2} dx\,dy \leq  C_r' \epsilon^{\delta}|\ln \epsilon|^d \eqno (1.5)$$
In the case where $K(x,y)g(z)$ is a smooth function of the form $\phi(x,y)$ and $\alpha = \beta = 0$, this $(\delta, d)$ will be the same as the  $(\delta, d)$ 
defined above except when the Hessian of $S(x,y)$ is nonvanishing. This relationship between sublevel set measures and oscillatory integral decay rate can be proven using resolution of
singularities and we again refer to [AGV] for details.

\noindent Another way to view $\delta$ is as follows. Define the meaure  $d\mu_{\alpha,\beta}$ by
$$d\mu_{\alpha,\beta} =  |S(x,y)|^{\alpha} (x^2 + y^2)^{\beta \over 2} \,dx
\,dy \eqno (1.6)$$
Then for sufficiently small $r > 0$, $\delta$ is also given by
$$\delta = \sup \{\eta: \int_{E_r} |S(x,y)|^{-\eta}\,\,d\mu_{\alpha,\beta} < \infty\} \eqno (1.7)$$
One can therefore refer to $\delta$ as the "critical integrability exponent" of $S(x,y)$ at $(0,0)$ with respect to the measure $d\mu_{\alpha,\beta}$.

\noindent Our main result is the following.

\noindent {\bf Theorem 1.1.} Let $o$ denote the order of the zero of $S(x,y)$ at the origin. 

\noindent {\bf a)} Suppose $\delta < {1 \over 3} + {1 \over 3o}$. There exists an $r > 0$ such that if $K(x,y)$ is  supported in $E_r$ then
$$|T(\lambda_1,\lambda_2,\lambda_3)| \leq C_{S,A}(1 + |\lambda_1|)^{-\delta} (\ln(1 + |\lambda_1|))^d \eqno (1.8a)$$
\noindent {\bf b}) Suppose $\delta > {1 \over 3} + {1 \over 3o}$. There exists an $r > 0$ such that if $K(x,y)$ is  supported in $E_r$ then
$$|T(\lambda_1,\lambda_2,\lambda_3)| \leq C_{S,A}  (1 + |\lambda_1|)^{-{1 \over 3} - {1 \over 3o}}\eqno (1.8b)$$
\noindent {\bf c)} Suppose $\delta = {1 \over 3} + {1 \over 3o}$. There exists an $r > 0$ such that if $K(x,y)$ is  supported in $E_r$ then  
$$|T(\lambda_1,\lambda_2,\lambda_3)| \leq C_{S,A}(1 + |\lambda_1|)^{-{1 \over 3} - {1 \over 3o}} (\ln(1 + |\lambda_1|))^{d+1} \eqno (1.8c)$$ 
Here the $A$ in $C_{S,A}$ is as in $(1.2a)-(1.2b)$.

In the case where $K(x,y)g(z)$ is a smooth function $\phi(x,y)$, by [V] equation $(1.8a)$ is best possible. In [V] it is also shown that in the case of smooth $\phi(x,y)$
one always has ${1 \over o} \leq \delta \leq {2 \over o}$. So $\delta$ will always be less than ${1 \over 3} + {1 \over 3o}$ if 
${2 \over o} < {1 \over 3} + {1 \over 3o}$, or equivalently if $o > 5$. Thus for smooth $\phi(x,y)$, whenever $o > 5$ case a) of  Theorem 1.1  is sharp.
The only situation where $o  = 5$ that Theorem 1.1 does not cover in the case of smooth $\phi(x,y)$  is when $\delta = {2 \over o}$, which only happens in special
situations and is easy to handle directly. For $o = 3$ and $4$, sometimes one is in case a) and other times one is in the nonsharp cases b) and
c). When $o = 2$, the statement that $\delta < {1 \over 3} + {1 \over 3o}$ reduces to the statement that $\delta < {1 \over 2}
= {1 \over o}$, which never occurs. 

It should be pointed out that in this smooth case, a weaker version of Theorem 1.1 was 
proven in Theorem 1.2 of [G1].

As we will see in section 6, when $\beta = 0$, the estimates provided by Theorem 1.1a) are sharp, and furthermore $(\delta,d) = (\alpha + \delta_0,d_0)$, where $(\delta_0,d_0)$ are the $(\delta,d)$ of the smooth case (with $\alpha = \beta = 0$).
We leave open the question of sharpness of the uniform estimates of Theorem 1.1 when $\beta \neq 0$. If the $\beta = 0$ 
case is any indication, part a) of Theorem 1.1 is more likely to be sharp than the other two parts.

In [V] it is shown that in the case of real analytic $S(x,y)$ and smooth $K(x,y)g(z)$ there is a nice description of $\delta$ and $d$ in terms of Newton polygons and adapted coordinates. This was generalized to smooth $S(x,y)$ in [IM1]. In the more general scenario of this paper, unfortunately such a description no longer holds, which is why we only use $(\delta,d)$ as defined here and do not delve  into Newton polygons and related matters.

\noindent {\bf 2. The resolution of singularities theorem and some consequences.}

Let $S(x,y)$ be any smooth function with $S(0,0) = 0$ such that the Taylor expansion of $S(x,y)$ does not vanish to infinite order at the origin. Let $o$ denote the order of the zero of $S(x,y)$ at $(0,0)$. After rotating coordinates if necessary, we may assume that the Taylor expansion $\sum_{\alpha ,\beta}s_{\alpha\beta}x^{\alpha}y^{\beta}$ of $S$ centered 
at the origin has a nonvanishing $s_{0\,o}x^o$ term and a nonvanishing $s_{o\,0}y^o$ term. In this paper, we will use the resolution of singularities theorem of [G1] (Theorem 2.1 of that paper). It proceeds as follows. First, one divides the
$xy$ plane into eight triangles by slicing the plane using  the $x$ and $y$ axes  and two lines through the origin, one of the form $y = mx$ for some $m > 0$ and one of the form $y = mx$ for some $m < 0$. One must ensure that these two lines are not ones  on 
which the function $S_0(x,y) = \sum_{\alpha + \beta = o} s_{\alpha\beta}x^{\alpha}y^{\beta}$ vanishes other than at the
origin. After reflecting about the $x$ and/or $y$ axes and/or the line $y = x$ if necessary, each of the triangles becomes of the form $T_b = \{(x,y) \in \R^2: x > 0,\,0 < y < bx\}$ (modulo an inconsequential boundary set of measure zero). Theorem 2.1 of
[G1] is then as follows.

\noindent {\bf Theorem 2.1.} Let  $T_b = \{(x,y) \in \R^2: x > 0,\,0 < y < bx\}$ be as above. Abusing notation slightly, use the notation $S(x,y)$ to denote the reflected function $S(\pm x,\pm y)$ or $S(\pm y, \pm x)$ corresponding to $T_b$.
 Then there is a $a > 0$ and a positive integer $N$ such that
if $F_a$ denotes  $\{(x,y) \in \R^2: 0 \leq x\leq a, \,0 \leq y \leq bx\}$, then one can write $F_a = \cup_{i=1}^n cl(D_i)$, such that for to each $i$ there is a $\phi_i(x)$ with $\phi_i(x^N)$ smooth and $\phi_i(0) = 0$ such that after a coordinate change of the form $\eta_i(x,y) = (x, \pm y + \phi_i(x))$, the set $D_i$ becomes a set $D_i'$ on which the function $S \circ \eta_i(x,y)$ approximately becomes a monomial $d_i x^{\alpha_i}y^{\beta_i}$, $\alpha_i$ a nonnegative rational number and $\beta_i$ a nonnegative integer as follows.

\noindent {\bf a)} $D_i' = \{(x,y): 0 < x < a, \, g_i(x) < y < G_i(x)\}$, where $g_i(x^N)$ and $G_i(x^N)$ are
smooth. If we expand $G_i(x) =  H_i x^{M_i} + ...$, then $M_i \geq 1$ and $H_i > 0$, and consists of a single term $H_ix^{M_i}$ when $\beta_i = 0$. 

\noindent {\bf b)} Suppose $\beta_i = 0$. Then $g_i(x) = 0$. Either $\phi_i(x) = k_ix$ for some $k_i$, 
or $\phi_i(x)$ is of the form $k_i x + l_ix^{s_i} + $ higher order terms (if any),  where $k_i,l_i\neq 0$ and $M_i \geq s_i > 1$. If there
are higher order terms then one has strict inequality $M_i > s_i$. In addition,  the set $D_i'$ can
be constructed such that for any predetermined $\eta > 0$ there is a $d_i \neq 0$ such that on $D_i'$, for all $0 \leq l \leq \alpha_i$ one has
$$ |\partial_x^l (S \circ \eta_i)(x,y) -  d_i\alpha_i(\alpha_i -1) ... (\alpha_i - l + 1)x^{\alpha_i - l}| < \eta |d_i|x^{\alpha_i-l} \eqno (2.1)$$
This $\eta$ can be chosen independent of all the exponents appearing in this theorem. 
Furthermore, if one Taylor expands $S \circ \eta_i(x,y)$ in powers of $x^{1 \over N}$ and $y$ as  $\sum_{\alpha,\beta}S_{\alpha,\beta}x^{\alpha}y^{\beta}$, then $\alpha_i \leq \alpha + M_i\beta$ for all $(\alpha,\beta)$ such that $S_{\alpha,\beta} \neq 0$, with equality holding for at least two
$(\alpha,\beta)$, one of which is $(\alpha_i,0)$ and another of which satisfies $\beta > 0$. 

\noindent {\bf c)} If $\beta_i > 0$, then either $g_i(x)$ is identically zero or $g_i(x)$ 
can be expanded as $h_ix^{m_i} + ...$ where $h_i > 0$ and $m_i > M_i$. In addition, one may write $S = S_1^i + S_2^i$ as follows. $S_2^i \circ \eta_i(x,y)$ has a zero of infinite order at $(0,0)$ and
is identically zero if $S$ is real-analytic. $S_1^i \circ \eta_i(x^N,y)$ is smooth and there exists a $d_i \neq 0$ such that for any predetermined $\eta > 0$ the $D_i'$ can
be constructed such that on $D_i'$, for any $0 \leq l \leq \alpha_i$ and any $0 \leq  m \leq \beta_i$ one has
$$|\partial_x^l\partial_y^m(S_1^i \circ \eta_i)(x,y) -  \alpha_i(\alpha_i - 1)....(\alpha_i - l + 1)\beta_i(\beta_i - 1)...(\beta_i - m + 1)
d_ix^{\alpha_i - l}y^{\beta_i - m}| $$
$$\leq \eta |d_i| x^{\alpha_i - l}y^{\beta_i - m} \eqno (2.2)$$

\noindent The next lemma shows that the $(\delta,d)$ of Theorem 1.1 is well-defined.

\noindent {\bf Lemma 2.2.} Suppose that $S(x,y)$ is real-analytic on a neighborhood of the origin with $S(0,0) = 0$, and let $E_r$ 
denote the disk $\{(x,y): x^2 + y^2 < r^2\}$. Let $\alpha$
and $\beta$ be real numbers such that $|S(x,y)|^{\alpha}(x^2 + y^2)^{\beta \over 2}$ is integrable on a neighborhood of 
the origin. Then there is a $\delta > 0$ and an integer $d = 0$ or $1$ such that if $r$ is sufficiently small then there are 
constants $C$ and $C'$ depending on $\alpha$, $\beta$, $S(x,y)$, and $r$ such that for $0 < \epsilon < {1 \over 2}$ one has
$$ C \epsilon^{\delta}|\ln \epsilon|^d  \leq \int_{\{(x,y) \in E_r :  |S(x,y)| < \epsilon\}}|S(x,y)|^{ \alpha} (x^2 + y^2)^{\beta \over 2} dx\,dy \leq  C'  \epsilon^{\delta}|\ln \epsilon|^d \eqno (2.3)$$
\noindent {\bf Proof.} Let $D_i$ and $D_i'$ be the domains coming from applying Theorem 2.1 to $S(x,y)$. Then it suffices to 
show $(2.3)$ holds with $E_r$ replaced by $E_r \cap D_i$ and then the result follows from addition. If one does the coordinate
change $\phi_i$ of Theorem 2.1, on the new domain $D_i'$ one has that $|S \circ \eta_i (x,y)|$ is within a constant factor of $x^{\alpha_i}
y^{\beta_i}$ and that $(x^2 + y^2)^{\beta \over 2}$ is within a constant factor of $x^{\beta}$.  Thus there are constants $c$ and $c'$
such that 
$$\int_{{\{(x,y) \in D_r':
0 < x < c r,\,x^{\alpha_i}y^{\beta_i} < c\epsilon \}}}x^{\alpha\alpha_i + \beta}y^{\alpha \beta_i} < \int_{\{(x,y) \in E_r :  |S(x,y)| < \epsilon\}}|S(x,y)|^{ \alpha} (x^2 + y^2)^{\beta \over 2} $$
$$ < \int_{{\{(x,y) \in D_r':
0 < x < c'r,\,x^{\alpha_i}y^{\beta_i} < c'\epsilon \}}}x^{\alpha\alpha_i + \beta}y^{\alpha \beta_i} \eqno (2.4)$$
Recall that the upper boundary of $D_i'$ is of the form $H_i x^{M_i} + ...$ and the lower boundary is of the form $h_ix^{m_i} + ...$
for $m_i > M_i$ or is the $x$-axis. In the former case, we define $F_{r,c}$ and $G_{r,c'}$ by
$$F_{r,c} = \{(x,y):
0 < x < c r,\,2h_ix^{m_i} < y  <  {H_i \over 2} x^{M_i} ,\,x^{\alpha_i}y^{\beta_i} < c\epsilon \} \eqno (2.5a)$$
$$G_{r,c'} = {\{(x,y): 0 < x < c'r,{h_i \over 2}x^{m_i} < y  <  2H_i} x^{M_i},\,x^{\alpha_i}y^{\beta_i} < c'\epsilon \} \eqno (2.5b)$$
So if $r$ is sufficiently small there are positive constants $c$ and $c'$ such that one has 
$$\int_{F_{r,c}}x^{\alpha\alpha_i + \beta}y^{\alpha \beta_i} < \int_{\{(x,y) \in E_r :  |S(x,y)| < \epsilon\}}|S(x,y)|^{ \alpha} (x^2 + y^2)^{\beta \over 2} < \int_{G_{r,c'}}x^{\alpha\alpha_i + \beta}y^{\alpha \beta_i} \eqno (2.6)$$
One can directly compute the integrals on the left and right of $(2.5)$ and one obtains
an expression of the form
$C\epsilon^{\delta}|\ln \epsilon|^d + o(\epsilon^{\delta}|\ln \epsilon|^d)$ for both of them for some $\delta > 0$ and $d = 0$ or $1$. These are the needed estimates for a $D_i'$ whose lower boundary is not the $x$-axis. If the lower boundary of $D_i'$ is
the $x$-axis, we take $h_i = 0$ in the above and we get the needed estimates once again. This completes the proof of Lemma 2.2.

\noindent {\bf 3. Proof of Theorem 1.1 when $\beta_i > 0$.} 

In this and future sections we will make frequent use of the following classical Van der Corput Lemma (see p 334 of [S]):

\noindent {\bf Lemma 3.1.} Suppose $P(x)$ is a $C^k$ function on the interval $[a,b]$ with $|P^{(k)}(x)| > M$ on $[a,b]$ for
some $M > 0$. Let $\psi(x)$ be $C^1$ on $[a,b]$. If $k \geq 2$ there is a constant $c_k$ depending only on $k$ such that
$$\bigg|\int_a^b e^{iP(x)}\psi(x)\,dx\bigg| \leq c_kM^{-{1 \over k}}\bigg(|\psi(b)| + \int_a^b |\psi'(x)|\,dx\bigg)$$
If $k =1$, the same is true if we add the conditions that $P(x)$ is $C^2$ and that $P'(x)$ is monotonic on $[a,b]$. 

We also will make use of the following variation of the classical Van der Corput Lemma that holds for mixed partial derivatives.

\noindent {\bf Lemma 3.2.} Let $I_1$ and $I_2$ be closed intervals of lengths $l_1$ and $l_2$ respectively, and for some 
strictly monotone functions $f_1(x)$ and $f_2(x)$ on $I_1$ with $f_1(x) \leq f_2(x)$ let $R = \{(x,y) \in I_1 \times I_2: f_1(x) \leq y \leq  f_2(x)\}$ (Note $R$ might just be $I_1 \times I_2$). Suppose for some $k \geq 2$, $P(x,y)$ is a $C^k$ function on $R$ such that for each $(x,y) \in R$ one has
$$|\partial_{xy} P(x,y)| > M\,\,\,\,\,\,\,\,\,\,\,\,\,\,\,\,{\rm and }\,\,\,\,\,\,\,\,\,\,\,\,\,\,\,
\partial_y^k P(x,y) \neq 0 \eqno (3.1)$$
Further suppose that $\Psi(x,y)$ is a function on $R$ that is $C^1$ in the $y$ variable for fixed $x$,  such that 
$$ |\Psi(x,y)| < N \,\,\,\,\forall x,y\,\,\,\,\,\,\,\,\,\,\,\,\,\,\,{\rm and }\,\,\,\,\,\,\,\,\,\,\,\,\,\,\,\,\int_{\{y: (x,y) \in R\}} |\partial_y\Psi(x,y)|\,dy< N \,\,\,\,\forall x \eqno (3.2)$$
If  $R' \subset R$ such that the intersection of $R'$ with each vertical line is either empty or is a set of at most $l$ intervals, then 
the following estimate holds. 
$$\bigg|\int_{R'} e^{i P(x,y)}\Psi(x,y)\,dx\,dy\bigg| < C_{kl}  N \bigg({l_1l_2 \over M}\bigg)^{1 \over 2} \eqno (3.3)$$
\noindent {\bf Proof.}
 Write $\int_{R'} e^{i\lambda P(x,y)}\Psi(x,y)\,dx\,dy = I_1 + I_2$, where 
$$ I_1 = \int_{\{(x,y) \in R': |\partial_y P(x,y)| < ({Ml_1 \over l_2})^{1 \over 2}\}}e^{iP(x,y)}\Psi(x,y)\,dx\,dy \eqno (3.4a)$$
$$ I_2 = \int_{\{(x,y) \in R': |\partial_y P(x,y)| > ({Ml_1 \over l_2})^{1 \over 2}\}}e^{iP(x,y)}\Psi(x,y)\,dx\,dy \eqno (3.4b)$$
We estimate  $ |I_1|$ simply by taking absolute values of the integrand and then integrating. Since $|\partial_x(\partial_y P(x,y))| > M$, for fixed $y$ the measure of the $x$ in $R$
for which $|\partial_yS(x,y)| < ({Ml_1 \over l_2})^{1 \over 2}$ is at most $({Ml_1 \over l_2})^{1 \over 2}\times {2 \over M} = 2({l_1 \over Ml_2})^{1 \over 2}$. Thus, using the left half of $(3.2)$, for fixed $y$ the $x$-integral in $(3.4a)$ is at most $N 2({l_1 \over Ml_2})^{1 \over 2}$. Integrating this in $y$ we see that
$$|I_1| <   2N \bigg({l_1l_2 \over M}\bigg)^{1 \over 2}\eqno (3.5)$$
These are the bounds we seek. 

We now move on to $I_2$. Note that since $\partial_y^k P(x,y) \neq 0$ on $R$, for fixed $x$ the set of $y \in R$ for which $|\partial_y P(x,y)| > ({Ml_1 \over l_2})^{1 \over 2}$ is the union of at most $k$ intervals. Thus for fixed $x$, the set of $y \in R'$ for which $|\partial_y P(x,y)| > ({Ml_1 \over l_2})^{1 \over 2}$ is at most $kl$ intervals. On each of these intervals
we use the Van der Corput Lemma 3.1 for first derivatives in the $y$ direction in conjunction with $(3.2)$, add up the results, and then integrate the result in $x$. Although $\partial_y P(x,y)$ is not necessarily monotone on each of the intervals and therefore Lemma 3.1 does not 
immediately apply, the fact that $\partial_y^k P(x,y) \neq 0$ on $R$ with $k \geq 2$ ensures that a given interval is the union of at most 
$k$ intervals on which $\partial_y P(x,y)$ is monotone and on which we can apply Lemma 3.1.

So using Lemma 3.1, we see that for given $x$ the absolute value of the $y$-integral in $(3.4b)$ is at most $C_{kl} N 
({l_2 \over Ml_1})^{1 \over 2}$. Integrating this in $x$ gives $C_{kl} N ({l_1 l_2 \over M})^{1 \over 2}$. So we have
$$|I_2| <   C_{kl}N \bigg({l_1l_2 \over M}\bigg)^{1 \over 2}\eqno (3.6)$$
Adding this to $(3.5)$ completes the proof of Lemma 3.2.

We now proceed to the proof of Theorem 1.1 for $\beta_i > 0$. We perform the resolution of singularities algorithm of Theorem 2.1 to $S(x,y)$, and correspondingly write $T(\lambda) = \sum_{i = 1}^n T_i(\lambda)$, where $T_i(\lambda)$ is given by
$$T_i(\lambda_1,\lambda_2,\lambda_3) = \int_{D_i} e^{i\lambda_1 S(x,y) + i\lambda_2 x + i\lambda_3 y}\,g(S(x,y))K(x,y)\,dx\,dy \eqno (3.7)$$
Shifting the $y$ variable by $\phi_i(x)$ as in Theorem 2.1, this becomes
$$T_i(\lambda_1,\lambda_2,\lambda_3) = \int_{D_i'}e^{i\lambda_1 S \circ \eta_i (x,y) + i\lambda_2 x \pm i\lambda_3 y + i\lambda_3\phi_i(x)}\,g(S \circ \eta_i (x,y))K \circ \eta_i (x,y)\,dx\,dy \eqno (3.8)$$
Without loss of generality, we will always  take $\pm i\lambda_3 y$ to be $i\lambda_3 y$. Note that by the form of $\phi_i(x)$ given by part d) of Theorem 2.1, $K_i(x,y) = K \circ \eta_i (x,y)$ satisfies 
$(1.2b)$. Writing  $S_i(x,y) = S \circ \eta_i (x,y)$ we have
$$T_i(\lambda_1,\lambda_2,\lambda_3) = \int_{D_i'} e^{i\lambda_1 S_i (x,y) + i\lambda_2 x + i\lambda_3 y + i\lambda_3\phi_i(x)}\,g(S_i (x,y))K_i (x,y)\,dx\,dy \eqno (3.9)$$
Let $(\alpha_i,\beta_i)$ be as in Theorem 2.1, so that $S_i(x,y)$ is within a bounded factor of $x^{\alpha_i}y^{\beta_i}$ on $D_i'$, with corresponding estimates for its derivatives. The analysis is broken up into three cases, when $\beta_i = 0$, when $\beta_i = 1$, and when $\beta_i > 1$, with the $\beta_i = 0$ case the hardest. We do the $\beta_i > 0$ cases in this section, and then do the $\beta_i = 0$ case in sections 4 and 5.

\noindent {\bf Case 1}. $\beta_i \geq 2$. 

We divide the domain of integration of $(3.9)$ dyadically in the $x$ and $y$ variables and correspondingly we write $T_i = \cup_{j,k} T_{ijk}$ where $T_{ijk}$ is given by
$$T_{ijk}(\lambda_1,\lambda_2,\lambda_3) = \int_{D_i' \cap [2^{-j-1},2^{-j}] \times [2^{-k-1},2^{-k}]}e^{i\lambda_1 S_i (x,y) + i\lambda_2 x + i\lambda_3 y + i\lambda_3\phi_i(x)}\,g(S_i (x,y))$$
$$\times K_i (x,y)\,dx\,dy \eqno (3.10)$$
Note that by Theorem 2.1 c) there is some constant $c$ depending only on $S(x,y)$ such that on the portion of $D_i'$ for which $ [2^{-j-1},2^{-j}] \times [2^{-k-1},2^{-k}]$  we have
$$|\partial_y^2\big(\lambda_1 S_i (x,y) + \lambda_2 x + \lambda_3 y + \lambda_3\phi_i(x)\big)|  = |\partial_y^2(\lambda_1 S_i (x,y))| $$
$$> c |\lambda_1|x^{\alpha_i}y^{\beta_i - 2} \eqno (3.11)$$
We will now use $(3.11)$ and apply the Van der Corput lemma, Lemma 3.1, in the $y$ direction in $(3.10)$. For this we need to bound the $y$ derivatives of $g(S_i(x,y))$ and 
$K_i(x,y)$. As mentioned above, $(1.2b)$ holds for $K_i(x,y)$ in place of $K(x,y)$, so we have
$$|\partial_y K_i(x,y)| \leq C(x^2 + y^2)^{{\beta - 1\over 2}}$$
$$\leq C'{1 \over y} (x^2 + y^2)^{{\beta \over 2}} \eqno (3.12)$$
For $g(S_i(x,y))$, note that we have 
$$\partial_y (g(S_i(x,y)))= g'(S_i(x,y)) \partial_y S_i(x,y) \eqno (3.13)$$
So by $(1.2a)$ and $(2.2)$ we have
$$|\partial_y (g(S_i(x,y)))| \leq C|S_i(x,y)|^{\alpha -1}x^{\alpha_i}y^{\beta_i - 1}$$
$$\leq C'{1 \over y} |S_i(x,y)|^{\alpha} \eqno (3.14)$$ 
Taking $(3.12)$ and $(3.14)$ together, we have that the factor $g(S_i (x,y))K_i (x,y)$ satisfies
$$|\partial_y [g(S_i (x,y))K_i (x,y)]| \leq C {1 \over y} |S_i(x,y)|^{\alpha}(x^2 + y^2)^{{\beta \over 2}} $$
On the other hand, by $(1.2a)-(1.2b)$ one has
$$|g(S_i (x,y))K_i (x,y)| \leq C |S_i(x,y)|^{\alpha}(x^2 + y^2)^{{\beta \over 2}}$$
On the support of the integrand of $(3.9)$, we have $0 < y < Cx$, $x \sim 2^{-j}$, and $y \sim 2^{-k}$, so the last two
equations can be rewritten as
$$ |\partial_y [g(S_i (x,y))K_i (x,y)]| \leq C (2^k)(2^{-j\alpha_i\alpha - k\beta_i\alpha})(2^{-j\beta})  \eqno (3.15a)$$
$$ |g(S_i (x,y))K_i (x,y)| \leq C (2^{-j\alpha_i\alpha - k\beta_i\alpha})(2^{-j\beta})  \eqno (3.15b)$$
We apply Lemma 3.1 in the $y$ direction in $(3.10)$, using  $(3.11)$, $(3.15a)$ and $(3.15b)$. We get that for fixed $x$, the $y$ integral is bounded by
$$C |\lambda_1|^{-{1 \over 2}}(2^{-j\alpha_i\alpha - k\beta_i\alpha})(2^{-j\beta})( 2^{{j\alpha_i + k\beta_i \over 2} - k}) \eqno (3.16b)$$
Integrating this in $x$ we obtain
$$|T_{ijk}(\lambda)| \leq C |\lambda_1|^{-{1 \over 2}} (2^{-j\alpha_i\alpha - k\beta_i\alpha})(2^{-j\beta})( 2^{{j\alpha_i + k\beta_i \over 2} - j - k}) \eqno (3.17)$$
Equation $(3.17)$ implies
$$|T_{ijk}(\lambda)| \leq C \int_{ [2^{-j-1},2^{-j}] \times [2^{-k-1},2^{-k}] }\,\,\,\,|\lambda_1|^{-{1 \over 2}}|x^{\alpha_i}y^{\beta_i}|^{\alpha}(x^2 + y^2)^{\beta}|x^{\alpha_i}y^{\beta_i}|^{-{1 \over 2}} \eqno (3.18)$$
Recalling the definition $(1.6)$ of the measure $d\mu_{\alpha,\beta}$, this is the same as
$$|T_{ijk}(\lambda)| \leq C\int_{ [2^{-j-1},2^{-j}] \times [2^{-k-1},2^{-k}] }\,\,\,\,|\lambda_1|^{-{1 \over 2}}|x^{\alpha_i}y^{\beta_i}|^{-{1 \over 2}} \,d\mu_{\alpha,\beta} \eqno (3.19)$$
By simply taking absolute values of the integrand in $(3.10)$ and integrating, in view of $(1.2a)-(1.2b)$ we have
$$|T_{ijk}(\lambda)| \leq A^2\int_{ [2^{-j-1},2^{-j}] \times [2^{-k-1},2^{-k}] }\,\,\,\,1 \, d\mu_{\alpha,\beta} \eqno (3.20)$$
Combining $(3.19)$ and $(3.20)$ one then has
$$|T_{ijk}(\lambda)| \leq C \int_{ [2^{-j-1},2^{-j}] \times [2^{-k-1},2^{-k}] }\,\,\,\,\min(1, |\lambda_1 x^{\alpha_i}y^{\beta_i}|^{-{1 \over 2}})\,d\mu_{\alpha,\beta} \eqno (3.21)$$
Adding $(3.21)$ over all $j$ and $k$, and using the shape of $D_i'$ given by Theorem 2.1, we obtain
$$|T_i(\lambda)| \leq C  \int_{D_i'} \min(1, |\lambda_1 x^{\alpha_i}y^{\beta_i}|^{-{1 \over 2}})\, d\mu_{\alpha,\beta} \eqno (3.22)$$
Since $S_i(x,y) \sim x^{\alpha_i}y^{\beta_i}$ on $D_i'$, $(3.22)$ implies 
$$|T_i(\lambda)| \leq C  \int_{D_i'} \min(1, |\lambda_1 S_i(x,y)|^{-{1 \over 2}}) \,d\mu_{\alpha,\beta} \eqno (3.23)$$
Because $o \geq 2$, we have ${1 \over 3} + {1 \over 3o} \leq {1 \over 2}$ and therefore
$$|T_i(\lambda)| \leq C  \int_{D_i'} \min(1, |\lambda_1 S_i(x,y)|^{-{1 \over 3} - {1 \over 3o}})\, d\mu_{\alpha,\beta} \eqno (3.24)$$
$$ = \bigg(\mu_{\alpha,\beta}\big(\{(x,y) \in D_i': |S_i(x,y)| < {1 \over |\lambda_1|}\}\big) $$
$$+  {1 \over |\lambda_1|^{{1 \over 3}+ {1 \over 3o}}}
\int_{\{(x,y) \in D_i':\,\, |S_i(x,y)| \geq {1 \over |\lambda_1|}\}} {1 \over |S_i(x,y) |^{{1 \over 3}+ {1 \over 3o}}}\,d\mu_{\alpha,\beta}\bigg) \eqno (3.25)$$
By the definition $(1.5)$ of $(\delta, d)$, the first term of $(3.25)$ is at most $C|\lambda_1|^{-\delta}
(\ln|\lambda_1|)^d$, which is at least as good as the estimates of Theorem 1.1 in all three cases. As for the second term, by the characterization of integrals in terms of distribution functions (applied to ${1 \over |S_i(x,y)|})$ we have
$\int_{\{(x,y) \in D_i':\,\, |S_i(x,y)| \geq {1 \over |\lambda_1|}\}} {1 \over |S_i(x,y) |^{{1 \over 3}+ {1 \over 3o}}}\,d\mu_{\alpha,\beta}$ is equal to
$$ \int_{1 \over |\lambda_1|}^{\infty}\bigg({1 \over 3}+ {1 \over 3o}\bigg)t^{-{4 \over 3} - {1 \over 3o}}\mu_{\alpha,\beta}\big(\{(x,y) \in D_i': {1 \over |\lambda_1|} < |S_i(x,y)| < t \}\big)\,dt\eqno (3.26a)$$
One can replace the upper bound of $\infty$ in $(3.26a)$ by just ${1 \over 2}$ as the difference results in a contribution of $C|\lambda_1|^{-{1 \over 3} + {1 \over 3o}}$ to the second term of $(3.25)$, which is always at least as good as the desired bound. Thus we will bound
$$ \int_{1 \over |\lambda_1|}^{1 \over 2} \bigg({1 \over 3}+ {1 \over 3o}\bigg)t^{-{4 \over 3} - {1 \over 3o}}\mu_{\alpha,\beta}\big(\{(x,y) \in D_i': {1 \over |\lambda_1|} < |S_i(x,y)| < t \}\big)\,dt\eqno (3.26b)$$
By $(1.5)$ the expression in $(3.26b)$ is bounded by
$$C \int_{1 \over |\lambda_1|}^{1 \over 2} \bigg({1 \over 3}+ {1 \over 3o}\bigg) t^{-{4 \over 3} - {1 \over 3o}}t^{\delta}(\ln t)^d\,dt\eqno (3.27)$$
If $\delta < {1 \over 3}+ {1 \over 3o}$, $(3.27)$ becomes $C'|\lambda_1|^{ {1 \over 3}+ {1 \over 3o} - \delta}\ln |\lambda_1|^d$
plus a smaller error term. Thus the
second term in $(3.25)$ is bounded by ${1 \over |\lambda_1|^{{1 \over 3}+ {1 \over 3o}}}$ times this, or $C|\lambda_1|^{- \delta}\ln |\lambda_1|^d$. Adding together with the first term of $(3.25)$, we see that
$$|T_i(\lambda)| \leq C|\lambda_1|^{- \delta}\ln |\lambda_1|^d \eqno (3.28)$$
Since by just taking absolute values of the integrand and integrating one has $|T_i(\lambda)|$ is bounded by a constant, one can 
also say that 
$$|T_i(\lambda)| \leq C(1 + |\lambda_1|)^{- \delta}\ln(1 +  |\lambda_1|)^d \eqno (3.29)$$
This gives the estimate $(1.8a)$ required by Theorem 1.1 for the situation where $\delta < {1 \over 3}+ {1 \over 3o}$.
Suppose now $\delta > {1 \over 3}+ {1 \over 3o}$. Then the expression $\mu_{\alpha,\beta}\big(\{(x,y) \in D_i': {1 \over |\lambda_1|} < |S_i(x,y)| < t \}\big)$ in $(3.26b)$ is bounded by $Ct^{\delta}|\ln t|^d$. So $(3.26b)$ is bounded by
$$ C\int_{1 \over |\lambda_1|}^{1 \over 2} t^{\delta -{4 \over 3} - {1 \over 3o}}|\ln t|^d
\,dt \eqno (3.30)$$
Since $\delta > {1 \over 3} + {1 \over 3o}$, $(3.30)$ is bounded by a constant. Hence the 
 the second term in $(3.25)$ is bounded by $C|\lambda_1|^{-{1 \over 3}- {1 \over 3o}}$, so for the $\delta < {1 \over 3}+ {1 \over 3o}$ case we get the estimate
$$|T_i(\lambda)| \leq C(1 + |\lambda_1|)^{-{1 \over 3}- {1 \over 3o}} \eqno (3.31)$$
This gives $(1.8b)$. Lastly, if $\delta = {1 \over 3} + {1 \over 3o}$, $(3.30)$ is bounded by a constant
times $(\ln|\lambda_1|)^{d+1}$, so putting this back into the second term of $(3.25)$ we now get
$$|T_i(\lambda)| \leq C(1 + |\lambda_1|)^{-{1 \over 3}- {1 \over 3o}}\ln(1 + |\lambda_1|)^{d+1}$$
This gives $(1.8c)$ and we are done with the proof of Theorem 1.1 for when $\beta \geq 2$.

\noindent {\bf Case 2.} $\beta_i = 1$. We once again write
$T_i = \sum_{j,k} T_{ijk}$ by dyadically decomposing in the $x$ and $y$ variables. So we have
$$ T_{ijk}(\lambda) = \int_{D_i' \cap [2^{-j-1},2^{-j}] \times [2^{-k-1},2^{-k}]} e^{i\lambda_1 S_i (x,y) + i\lambda_2 x + i\lambda_3 y + i\lambda_3\phi_i(x)}g(S_i (x,y))K_i (x,y)\, dx\,dy \eqno (3.32)$$
We apply Lemma 3.2 to the integral in $(3.32)$. If $P(x,y)$ denotes the phase function in $(3.32)$, we have
$$|\partial_{xy}P(x,y)| = |\lambda_1 \partial_{xy}S_i(x,y)|$$
$$> C|\lambda_1 x|^{\alpha_i - 1} \eqno (3.33)$$
The last inequality follows from $(2.2)$. By $(1.2a)-(1.2b)$ and the fact that $0 < y < Cx$ we have
$$|g(S_i (x,y))K_i (x,y)| < C(x^{\alpha_i}y^{\beta_i})^{\alpha}(x^2 + y^2)^{\beta \over 2}$$
$$\leq C' (2^{-j\alpha_i\alpha - k\alpha})( 2^{-j\beta})  \eqno (3.34)$$
Exactly as in $(3.15a)$, we have
$$|\partial_y (g(S_i (x,y))K_i (x,y))| < C (2^k)(2^{-j\alpha_i\alpha - k\alpha})( 2^{-j\beta})  \eqno (3.35)$$
Thus, as needed for Lemma 3.2, we have an estimate for $|\partial_y (g(S_i (x,y))K_i (x,y))|$ that is $2^k$ times the
estimate for $|g(S_i (x,y))K_i (x,y)|$ given in $(3.34)$. Applying Lemma 3.2 now, we get
$$|T_{ijk}(\lambda)| \leq C(2^{-j\alpha_i\alpha - k\alpha})( 2^{-j\beta}) {2^{-{j + k \over 2}} \over |\lambda_1|^{1 \over 2}
 2^{-j{\alpha_i} + j \over 2}} \eqno (3.36)$$
$$= C(2^{-j\alpha_i\alpha - k\alpha})( 2^{-j\beta}) {2^{-j - k} \over |\lambda_1|^{1 \over 2}
 2^{-j{\alpha_i} - k\over 2}} \eqno (3.37)$$
$$\leq C' \int_{[2^{-j-1},2^{-j}] \times [2^{-k-1},2^{-k}]} |\lambda_1|^{-{1 \over 2}}2^{{j\alpha_i + k \over 2}}d\mu_{\alpha,\beta} \eqno (3.38)$$
$$\leq C''  \int_{[2^{-j-1},2^{-j}] \times [2^{-k-1},2^{-k}]}  |\lambda_1|^{-{1 \over 2}}|x^{\alpha_i}y|^{-{1 \over 2}}d\mu_{\alpha,\beta} \eqno (3.39)$$
This is analogous to $(3.19)$ and is the estimate we seek. The argument from $(3.20)$ to $(3.31)$ now 
completes the proof for the case where $\beta_i = 1$.

\noindent {\bf 4. The $\beta_i = 0$ case away from the zeroes of $\lambda_1 \partial_{xx}S_i(x,y) + \lambda_3 \phi''(x)$.}

When we are away from the zeroes of  $\lambda_1 \partial_{xx}S_i(x,y) + \lambda_3 \phi''(x)$ in $(3.9)$, the argument resembles the argument when $\beta_i > 1$
except we apply the Van der Corput lemma in the $x$-direction instead of the $y$-direction. Write $\phi_i(x) = k_i x + \psi_i(x)$ where $\psi_i'(0) = 0$ and let $\lambda_4 = \lambda_2 + \lambda_3k_i$.
We divide dyadically in the $x$ variable, writing $T_i = \cup_j T_{ij}$, where 
$$T_{ij}(\lambda) = \int_{\{(x,y) \in D_i': \,x \in [2^{-j-1}, 2^{-j}]\}} e^{i\lambda_1 S_i (x,y) + \lambda_4x + i\lambda_3y  + i\lambda_3\psi_i(x)}\,g(S_i (x,y))K_i (x,y)\,dx\,dy \eqno (4.1)$$
Note that $\psi_i(x) = 0$ if $\phi_i(x)$ is linear, and $\psi_i(x) = l_i x^{s_i} + ...$ with $l_i \neq 0, 1 < s_i < \alpha_i$ otherwise.
 if $P_i(x,y)$ denotes the phase function in $(4.1)$, we have
$$\partial_{xx} P_i(x,y) = \lambda_1 \partial_{xx} S_i (x,y) + \lambda_3 \psi_i''(x) \eqno (4.2)$$
By Theorem 2.1b), $\partial_{xx} S_i (x,y)$ is equal to $d_i\alpha_i(\alpha_i - 1)x^{\alpha_i - 2}$ plus an error term less than $\eta x^{\alpha_i - 2}$ for an $\eta$ of our choice. Furthermore, when $\psi_i(x)$
is nonzero we have that  $\partial_{xx} \psi_i(x) = l_is_i(s_i - 1)x^{s_i - 2} + O(x^{s_i - 2 + \epsilon})$ for some $\epsilon > 0$. Because the two exponents $\alpha_i - 2$ and $s_i - 2$ are distinct, if $\eta$ is chosen appropriately then
for some constant $C$, for all but at most two values of $j$, on the support of the integrand of $(4.1)$ one has
$$ |\lambda_1 \partial_{xx} S_i (x,y) + \lambda_3 \partial_{xx} \psi_i(x)| > C(|\lambda_1 \partial_{xx} S_i (x,y)| + |\lambda_3 \partial_{xx} \psi_i(x)|) $$
$$ > C'|\lambda_1 \partial_{xx} S_i (x,y)| \eqno (4.3)$$
$$> C''|\lambda_1|  x^{\alpha_i-2} \eqno (4.4)$$
The latter inequality follows from  $(2.1)$. Note that if $\psi_i(x)$ is identically zero then $(4.4)$ holds for all $j$. Thus if $\psi(x)$ is not identically zero and $j$ is not one of these two exceptional values, or if $\psi(x)$ is identically zero and $j$ is anything, we may
argue as follows. We apply Lemma 3.1  in the $x$ direction in $(4.1)$, using $(4.4)$ on the phase and $(1.2a)-(1.2b)$ to bound
$g(S_i (x,y))K_i (x,y)$ and its $x$-derivative, analogously to as done in section 3. We then integrate the result in $y$, thereby gaining 
an additional factor of $C 2^{-j M_i}$, $M_i$ is as in Theorem 2.1, since $D_i'$ is of $y$-width comparable to $x^{M_i}$ for a given $x$. We get
$$|T_{ij}(\lambda)| < C_3(2^{-j\alpha_i\alpha})(2^{-j\beta}) |\lambda_1|^{-{1 \over 2}}( 2^{{j \alpha_i \over 2} - j})(2^{-M_i j})  \eqno (4.5)$$
Using the definition $(1.6)$ of $d\mu_{\alpha,\beta}$ and keeping in mind the area of the portion of $D_i'$ with $x \in [2^{-j-1}, 2^{-j}]$ is $\sim 2^{-jM_i - j}$, $(4.5)$ gives 
$$|T_{ij}(\lambda)| < C_4 \int_{\{(x,y) \in D_i': x \in [2^{-j-1},2^{-j}]\}}\,\,\,\,{1 \over |\lambda_1|^{1 \over 2}
2^{-{j\alpha_i \over 2 }}}\,d\mu_{\alpha,\beta} \eqno (4.6)$$
$$\leq  C_5 \int_{\{(x,y) \in D_i': x \in [2^{-j-1},2^{-j}]\}}\,\,\,\,{1 \over  |\lambda_1|^{1 \over 2}|S_i(x,y)|^{{1 \over 2}}} \,d\mu_{\alpha,\beta} \eqno (4.7)$$
Now we argue as in $(3.20)-(3.31)$ of Case 1 to achieve the desired estimates $(3.27)-(3.29)$.

Suppose now $\psi_i(x)$ is not identically zero and we are in one of the at most two exceptional $j$'s for which the above argument doesn't hold. If $N$ is any constant that depends only on $S(x,y)$, then by choosing the $\eta$ in Theorem 2.1b) according to $N$, one can cause $(4.4)$ to hold (with a different constant) outside of at most a vertical strip of width ${2^{-j} \over N}$. Namely, we define $x_0$ by the condition
$$\lambda_1 d_i\alpha_i(\alpha_i - 1)x_0^{\alpha_i - 2}+ \lambda_3  l_i s_i(s_i- 1)x_0^{s_i - 2} = 0 \eqno (4.8)$$
Then if $\eta$ is small enough, $(4.4)$ will hold outside of the set of $x$ where $|x - x_0| < {2^{-j} \over N}$. Then the the portion of $(4.1)$ outside  $|x - x_0| < {2^{-j} \over N}$
will once again satisfy the bounds of $(4.7)$ (with the constant depending on $N$). 
Thus in what follows, it suffices to bound $T_{ij}^N(\lambda)$, where $T_{ij}^N(\lambda)$ is given by
$$T_{ij}^N(\lambda) = \int_{D_i'}e^{i\lambda_1 S_i (x,y) + i\lambda_4x + i\lambda_3y  + i\lambda_3\psi_i(x)}\,g(S_i (x,y))K_i (x,y) \rho(2^jN(x - x_0))\,dx\,dy \eqno (4.9)$$
Here $\rho(x)$ is nonnegative, smooth, supported on $[-2,2]$, and equal to 1 on $[-1,1]$. The value of $N$ will be dictated by our arguments as we proceed, but will depend only on $S(x,y)$.
Note that $(4.8)$ gives an expression for $\lambda_3$ as a multiple of $\lambda_1$ as
$$\lambda_3 = -{d_i\alpha_i(\alpha_i - 1) \over l_is_i(s_i - 1)}x_0^{\alpha_i - s_i} \lambda_1 \eqno (4.10)$$

\noindent {\bf 5. The $\beta_i = 0$ case near the zeroes of $\lambda_1 \partial_{xx}S_i(x,y) + \lambda_3 \psi''(x)$.}

\noindent We start by writing the Taylor expansion of $S_i(x,y)$ in powers of $x^{1 \over N}$ and $y$ as 
$$S_i(x,y) = \sum_{\alpha, \beta} S_{\alpha \beta} x^{\alpha}y^{\beta} \eqno (5.1)$$
By Theorem 2.1, in the $\beta_i = 0$ situation, $D_i'$ is of the form $\{(x,y): 0 < x < a,\,\, 0 < y < H_ix^{M_i}\}$ for some
rational $M_i  \geq 1$. Thus it is natural to look at the function $S_i(x,x^{M_i}y)$ on the rectangle
$R_i = \{(x,y): 0 < x < a, \,\,0 < y < H_i\}$. On $R_i$, the Taylor expansion $(5.1)$ becomes
$$S_i(x,x^{M_i}y) = \sum_{\alpha, \beta} S_{\alpha \beta} x^{\alpha+M_i\beta }y^{\beta} \eqno (5.2)$$
In view of the form $(5.2)$, it makes sense to look at the (finitely many) terms of $(5.2)$ where $\alpha + M_i\beta$ takes its minimal value. By $(2.1)$, $S_i(x,y) \sim x^{\alpha_i}$ on $D_i'$, and taking $x \rightarrow 0$ in $(5.2)$ we see that $\alpha_i$ must be that minimal value. Thus if we let $p_i(x,y) = \sum_{\alpha + M_i\beta = \alpha_i} S_{\alpha \beta}x^{\alpha}
y^{\beta}$, then $(5.1)-(5.2)$ become
$$S_i(x,y) = p_i(x,y) + \sum_{\alpha + M_i\beta > \alpha_i} S_{\alpha \beta} x^{\alpha}y^{\beta} \eqno (5.3)$$
$$S_i(x,x^{M_i}y) = x^{\alpha_i}p_i(1,y) +  \sum_{\alpha + M_i\beta > \alpha_i} S_{\alpha \beta} x^{\alpha+M_i\beta }y^{\beta} \eqno (5.4)$$
By $(5.4)$, we have 
$$\lim_{x \rightarrow 0^+} {S_i(x,x^{M_i}y) \over x^{\alpha_i}} = p_i(1,y) \eqno (5.5)$$
Since by $(2.1)$ there are constants $C$ and $C'$ such that $Cx^{\alpha_i} < S_i(x,y) < C'x^{\alpha_i}$ on $D_i'$, 
$(5.5)$ implies that $p_i(1,y)$ has no zeroes for $y \in [0,H_i]$. Next, we look at $ \partial_y (S_i(x,x^{M_i}y))$, given by 
$$ \partial_y (S_i(x,x^{M_i}y)) = x^{\alpha_i} \partial_y p_i(1,y) + \sum_{\alpha + M_i\beta > \alpha_i}\beta S_{\alpha \beta} x^{\alpha+M_i\beta }y^{\beta - 1} \eqno (5.6)$$
By the last sentence of part b) of Theorem 2.1, $p_i(1,y)$ is nonconstant, so $\partial_y p_i(1,y)$ is not identically zero. 
Let $\alpha_i + \xi$ denote the minimal $\alpha + M_i\beta$ other than $\alpha_i$. So $\xi > 0$ and $(5.6)$ can be rewritten
as 
$$ \partial_y (S_i(x,x^{M_i}y)) = x^{\alpha_i} \partial_y p_i(1,y) +  x^{\alpha_i + \xi} q(x,y) \eqno (5.7)$$
Here $q(x,y)$ is a real-analytic function of $x^{1 \over N}$ and $y$. Next, we write $T_{ij}^N(\lambda)$ in the new coordinates:
$$T_{ij}^N(\lambda) = \int_{[2^{-j-1}, 2^{-j}] \times [0,H_i]}e^{i\lambda_1 S_i (x,x^{M_i}y) + i\lambda_4x + i\lambda_3x^{M_i} y  + i\lambda_3\psi_i(x)}\,g(S_i (x,x^{M_i}y))K_i (x,x^{M_i}y)$$
$$\times  \rho(2^jN(x - x_0))x^{M_i}\,dx\,dy \eqno (5.8)$$
If $\partial_y p_i(1,y)$ has any (real) roots in $[0,H_i]$, we enumerate them as $r_1,...,r_K$ and let $I_k$ denote the interval $[r_{ik} - \epsilon_0, r_{ik} + \epsilon_0]$, where $\epsilon_0$ denotes a small constant to be determined by our arguments. We write $T_{ij}^N(\lambda) = T_{ij}^{N,1}(\lambda) + T_{ij}^{N,2}(\lambda)$, where 
$$T_{ij}^{N,1}(\lambda) =  \int_{[2^{-j-1}, 2^{-j}] \times ([0,H_i] - \cup_k I_k)}e^{i\lambda_1 S_i (x,x^{M_i}y) + i\lambda_4x + i\lambda_3x^{M_i} y  + i\lambda_3\psi_i(x)}\,g(S_i (x,x^{M_i}y))$$
$$\times  K_i (x,x^{M_i}y)\rho(2^jN(x - x_0))x^{M_i}\,dx\,dy \eqno (5.9a)$$
$$T_{ij}^{N,2}(\lambda) =  \int_{[2^{-j-1}, 2^{-j}] \times ([0,H_i]  \cap (\cup_k I_k))}e^{i\lambda_1 S_i (x,x^{M_i}y) + i\lambda_4x + i\lambda_3x^{M_i} y  + i\lambda_3\psi_i(x)}\,g(S_i (x,x^{M_i}y))$$
$$\times  K_i (x,x^{M_i}y)\rho(2^jN(x - x_0))x^{M_i}\,dx\,dy \eqno (5.9b)$$
In the situation where $\partial_y p_i(1,y)$ has no roots in $[0,H_i]$, we just set $T_{ij}^{N,1}(\lambda) = T_{ij}^N (\lambda)$ and 
$T_{ij}^{N,2}(\lambda) = 0$. 

\noindent {\bf Estimates when $y$ is away from the zeroes of $\partial_y p_i(1,y)$}.

We now bound $T_{ij}^{N,1}(\lambda)$ through an application of Lemma 3.1 or Lemma 3.2. By Theorem 2.1, we always have
$M_i \geq s_i$. For now we assume $M_i > s_i$ and at the end of the argument we will describe the modifications needed for the $M_i = s_i$ situation. Let
$P(x,y)$ denote the phase function of $(5.9a)-(5.9b)$. Then using $(5.7)$ we have
$$\partial_{xy} P(x,y) = \lambda_1 \alpha_i x^{\alpha_i - 1}\partial_y p_i(1,y) + \lambda_1 \partial_x (x^{\alpha_i  + \xi}q(x,y)) + \lambda_3M_i x^{M_i - 1} \eqno (5.10)$$
Let $\epsilon_1 > 0$ be such that $|\partial_y p_i(1,y)| > \epsilon_1$ on $[0,H_i] - \cup_i I_i$. Since $\xi > 0$, if $x$ is small enough, which we may assume, then the $\lambda_1 \partial_x (x^{\alpha_i  + \xi}q(x,y))$ term in $(5.10)$ will be less than 
${\epsilon_1 \over 3}|\lambda_1\alpha_ix^{\alpha_i  - 1}|$ in absolute value. Furthermore, by $(4.10)$ one has $|\lambda_3M_i x^{M_i - 1}| <  C|\lambda_1 x^{\alpha_i - s_i + M_i - 1}|$. Since $\alpha_i > M_i > s_i$ in situation at hand, when $x$
 is sufficiently small the $\lambda_3M_i x^{M_i - 1}$ term is also of absolute value less than ${\epsilon_1 \over 3}|\lambda_1 \alpha_i x^{\alpha_i  - 1}|$. Thus when $x$ is sufficiently small, $(5.10)$ gives 
$$|\partial_{xy}P(x,y)| > {\epsilon_1 \over 3}|\lambda_1 \alpha_i x^{\alpha_i - 1}| \eqno (5.11a)$$
Since $x \sim 2^{-j}$ on the domain of integration of $(5.9a)$, we can rewrite $(5.11a)$ as
$$|\partial_{xy}P(x,y)| > C  |\lambda_1| 2^{-j\alpha_i + j} \eqno (5.11b)$$
We now apply Lemma 3.2 as follows. The domain
of integration of $(5.9a)$ is the union of finitely many rectangles, on each of which we apply Lemma 3.2 using $(5.11b)$. 
By $(1.2a)-(1.2b)$ as in $(3.32)$, the function $A(x,y) = g(S_i (x,x^{M_i}y))K_i (x,x^{M_i}y)\rho(2^jN(x - x_0))x^{M_i}$ satisfies 
$$|A(x,y)| \leq C (2^{-j\alpha_i\alpha})( 2^{-j\beta}) 
(2^{-j M_i}) \eqno (5.12)$$
Taking a $y$ derivative of $g(S_i (x,x^{M_i}y))$ gives $ g'(S_i(x,x^{M_i}y))\partial_y(S_i(x,x^{M_i}y))$, so in view of $(1.2a)$ and the fact that  $S_i(x,x^{M_i}y) \leq Cx^{\alpha_i}$ we have
$$|\partial_y (g(S_i (x,x^{M_i}y)))| \leq C 2^{-j\alpha_i(\alpha - 1)}|\partial_y(S_i(x,x^{M_i}y))| \eqno (5.13a)$$
So by $(5.7)$ we have
$$|\partial_y \big(g(S_i (x,x^{M_i}y))\big)| \leq C 2^{-j\alpha_i\alpha} \eqno (5.13b)$$
Next, note that $\partial_y \big(K_i(x,x^{M_i}y)\big) = x^{M_i} (\partial_y K_i)(x,x^{M_i}y)$, so in view of $(1.2b)$ we have
$$|\partial_y \big(K_i(x,x^{M_i}y)\big)| \leq C2^{-jM_i - j\beta + j} \eqno (5.14a)$$
Since $M_i \geq 1$, in particular we have
$$|\partial_y \big(K_i(x,x^{M_i}y)\big)| \leq C2^{- j\beta} \eqno (5.14b)$$
Using $(5.13b)$ and $(5.14b)$, we see that
$$|\partial_y A(x,y)| \leq C (2^{-j\alpha_i\alpha})( 2^{-j\beta}) (2^{-j M_i}) \eqno (5.15)$$
We now apply Lemma 3.2 on each of the rectangles in $(5.9a)$, using $(5.11b)$ on the phase and $(5.12), (5.15)$ on 
$A(x,y)$. We then add the estimates over the various rectangles. The result is
$$|T_{ij}^{N,1}(\lambda)| \leq C (2^{-j\alpha_i\alpha})( 2^{-j\beta}) (2^{-j M_i}){2^{-{j \over 2}} \over |\lambda_1|^{{1 \over 2}}2^{{1 \over 2}(-j\alpha_i + j)}} \eqno (5.16)$$
$$= C(2^{-j\alpha_i\alpha})( 2^{-j\beta}) (2^{-jM_i - j}){1 \over |\lambda_1|^{1 \over 2}2^{-{j\alpha_i \over 2}}} \eqno (5.17a)$$
By simply taking absolute values in $(5.9a)$ and integrating, using $(5.12)$ one also has
$$|T_{ij}^{N,1}(\lambda)| \leq C(2^{-j\alpha_i\alpha})( 2^{-j\beta}) (2^{-jM_i - j}) \eqno (5.17b)$$
Combining $(5.17a)$ and $(5.17b)$ we see that
$$|T_{ij}^{N,1}(\lambda)| \leq C(2^{-j\alpha_i\alpha})( 2^{-j\beta}) (2^{-jM_i - j}) \min\bigg(1, {1 \over |\lambda_1|^{1 \over 2}2^{-{j\alpha_i \over 2}}}\bigg) \eqno (5.18a)$$
$$\leq C' \int_{\{(x,y) \in D_i':\,x \in [2^{-j-1},2^{-j}]\}}\min\bigg(1, {1 \over |\lambda_1|^{1 \over 2}|S_i(x,y)|^{1 \over 2}}\bigg)\,d\mu_{\alpha,\beta} \eqno (5.18b)$$
This is the estimate for $|T_{ij}^{N,1}(\lambda)|$ that we will need.

The above assumed that $M_i > s_i$. When $M_i = s_i$, the reason that the above argument doesn't always work is that in
$(5.10)$, the terms  $\lambda_1 \alpha_i x^{\alpha_i - 1}\partial_y p_i(1,y)$ and $\lambda_3M_i x^{M_i - 1}$ may cancel each other (although
the  $\lambda_1 \partial_x (x^{\alpha_i  + \xi}q(x,y))$ term is smaller than each of them). In this situation we use Lemma 3.1 for
first derivatives instead of Lemma 3.2 as above. Since $\alpha_i > M_i$, the corresponding terms $\lambda_1 x^{\alpha_i}\partial_y p_i(1,y)$ and $\lambda_3 x^{M_i}$  of $\partial_xS(x,y)$ will not cancel each other in the narrow region near $x_0$ in the support
of the integrand of $(5.9a)$, and the $\lambda_1 x^{\alpha_i  + \xi}q(x,y)$ term will be small compared to each of them. Thus we  may apply Lemma 3.1 in this fashion.

\noindent {\bf Estimates when $y$ is near a zero of $\partial_y p_i(1,y)$ of order 1}.

We now start bounding $|T_{ij}^{N,2}(\lambda)|$. We will bound the portion of the integral $(5.9b)$ over $[0,a] \times
([0,H_i] \cap I_k)$ for each $k$. Denote this integral by $J_{ijk}$. We first consider the case where $ \partial_y p_i(1,y)$ has a zero of order 1 at $y = r_{ik}$. In this case, by $(5.7)$, if $\epsilon_0$ and $a$ are sufficiently small, which we may assume, then on $[0,a] \times  ([0,H_i] \cap I_k)$ one has
$$|\partial_{yy} (S_i(x,x^{M_i}y)) | > C x^{\alpha_i} \eqno (5.19)$$
Thus we may apply Lemma 3.1 in the $y$ direction, using $(5.12)$ and $(5.15)$ on $A(x,y)$ and $(5.19)$ on the phase, and then integrate the result in $x$. We obtain 
$$|J_{ijk}| \leq C  (2^{-j\alpha_i\alpha})( 2^{-j\beta}) (2^{-j M_i}){1 \over |\lambda_1|^{-{1 \over 2}}2^{-{j\alpha_i \over 2}}}(2^{-j}) \eqno (5.20)$$
This is exactly the same as $(5.17a)$, so once again we get the bound given by $(5.18b)$.

\noindent {\bf A second resolution of singularities when $y$ is near a zero of $\partial_y p_i(1,y)$ of order greater than 1}.

When we bound the integrals $J_{ijk}$ for the $r_{ik}$ at which  $\partial_y p_i(1,y)$ has a zero of order greater than 
1, the argument is more involved and uses a second application of resolution of singularities, this time to the function $\partial_{yy} (S_i(x,x^{M_i}y))$ on 
a square centered at $(0,r_{ik})$. To be precise, we shift the $y$-coordinate by $r_{ik}$ and then apply resolution of singularities to 
  $\partial_{yy} (S_i(x,x^{M_i}(y + r_{ik})))$ on a sufficiently small square $U_{ik}$  centered at the origin. The version of Theorem 2.1 is not exactly the one we need here. Instead we use the 
(very closely related) Theorem 3.1 of [G2], which says that if $U_{ik}$ is a sufficiently small enough square centered at $(0,0)$, 
for any $\eta > 0$ the portion of the square where where $|y| < |x|^{\eta}$ can be subdivided as in Theorem 2.1. (The proof
of Theorem 3.1 of  [G2] works the same way for an analytic function of $x^{1 \over N}$ and $y$ as it does for 
an analytic function of $x$ and $y$).  Although 
only a weaker version of $(2.1)$ is given in [G2] in this set-up, $(2.1)$ still holds here for exactly the same reason it holds in the
setup of Theorem 2.1; we omit the details for brevity. Also, one has to replace the conclusion $M_i \geq 1$ in part a) of
Theorem 2.1 with $M_i \geq \eta$ in this situation. 

In the coordinates of this application of Theorem 3.1 of [G2], the set $[0,a] \times ([0,H_i] \cap I_k)$ becomes of the form
$[0,a] \times [-\epsilon_0, \epsilon_0]$ if the root $r_{ik}$ is not an endpoint of the interval $[0,H_i]$, and is of the form 
$[0,a] \times [0, \epsilon_0]$ or $[0,a] \times [-\epsilon_0, 0]$ if it is. We select $\epsilon_0$ to be small enough to be the 
radius of a square on which the above resolution of singularities algorithm holds for each $r_{ik}$, and then set $a = \epsilon_0$.
 (We can shrink $a$ as much as we like as this only shrinks the neighborhood of the origin on which Theorem 1.1 holds.)
Then in the new coordinates the set $[0,a] \times ([0,H_i] \cap I_k)$ is either $[0,\epsilon_0] \times  [-\epsilon_0, \epsilon_0]$,
 $[0,\epsilon_0] \times  [0, \epsilon_0]$, or  $[0,\epsilon_0] \times  [-\epsilon_0, 0]$. Since the resolution of singularities 
algorithm of Theorem 3.1 of [G2] starts by dividing into 4 squares via the $x$ and $y$ axes and then does further subdivisions 
afterwards, the resolution of singularities procedure restricted to $[0,a] \times ([0,H_i] \cap I_k)$ will simply result in a 
subset of the set of domains given by the overall procedure. 

We let $D_{ikl}$ and $D_{ikl}'$  denote the domains for this second resolution of singularities that are analogous to the $D_i$ and $D_i'$ in Theorem 2.1,
and we let $\phi_{ikl}$ be the analogues of the coordinate changes $\phi_i$.
We denote the transformed $\partial_{yy} (S_i(x,x^{M_i}(y + r_{ik})))$ in the new coordinates by $Q_{ikl}(x,y)$. Let
$\alpha_{ikl}$ and $\beta_{ikl}$ denote the analogues of the exponents $\alpha_i$ and $\beta_i$, so that
$\partial_{yy} (S_i(x,x^{M_i}(\pm y + r_{ik} + \phi_{ikl}(x)))) \sim x^{\alpha_{ikl}}y^{\beta_{ikl}}$ on  $D_{ikl}'$, with corresponding estimates on 
its derivatives. Because  $x^{\alpha_i}$ divides  $\partial_{yy} (S_i(x,x^{M_i} y))$ by $(5.7)$ and because the coordinate change is of the 
form $(x,y) \rightarrow (x, \pm y + \phi_{ikl}(x))$, $x^{\alpha_i}$  also divides $Q_{ikl}(x,y)$, 
so that $\alpha_{ikl} \geq \alpha_i$ and we may write $Q_{ikl}(x,y) = x^{\alpha_i}\tilde{Q}_{ikl}(x,y)$ for some
function $\tilde{Q}_{ikl}(x,y)$ which is a real-analytic function of $x^{1 \over N'}$ and $y$ for some positive integer $N'$. Similarly, in the
new coordinates  $S_i(x,x^{M_i}(y + r_{ik}))$  can be written in the form $x^{\alpha_i}S_{ikl}(x,y)$. Note that due to the form of the
coordinate change we have $\partial_{yy} (x^{\alpha_i}S_{ikl}(x,y)) = x^{\alpha_i}\tilde{Q}_{ikl}(x,y)$,
or in other words, $\partial_{yy} (S_{ikl}(x,y)) = \tilde{Q}_{ikl}(x,y)$.

We also incorporate the $(x,y)$ for which $|y| > x^{\eta}$ into the collection of $D_{ikl}'$. We do this by simply by letting one $D_{ikl}'$ be $\{(x,y): 0 < x < \epsilon_0, \,x^{\eta} < y < \epsilon_0\}$
and another be  $\{(x,y): 0 < x < \epsilon_0,  -\epsilon_0 <  y <  -x^{\eta}  \}$, and then let the function $\phi_{ikl}(x)$ just be $0$. In this case, assuming $\eta$ is chosen 
sufficiently small, by $(5.7)$ there will again be an $\alpha_{ikl}$ and $\beta_{ikl}$ such that $Q_{ikl}(x,y) \sim x^{\alpha_{ikl}}y^{\beta_{ikl}}$ on $D_{ikl}'$, if we set $\alpha_{ikl} = \alpha_i$.

\noindent We have the following fact about  $x^{\alpha_{ikl}}y^{\beta_{ikl}}$ .

\noindent {\bf Lemma 5.1.} There is a constant $C$ such that $x^{\alpha_{ikl}}y^{\beta_{ikl}} > Cx^{\alpha_i}y^{o - 2}$ on $D_{ikl}'$ for each $(i,k,l)$.

\noindent {\bf Proof.} We start by noting that there is a
constant $c$ such that  $|\partial_y^o S(x,y)| > c$ and $|\partial_x^o S(x,y)|  > c$ in the original (rotated) coordinates on a neighborhood of the origin, and $S_i(x,y)$ is either $S(x,y)$ after a coordinate change of the form $(x,y) \rightarrow (\pm x, \pm y + \phi_i(x))$, or is $ S(x,y)$ after doing the coordinate change $(x,y) \rightarrow (y,x)$ and then doing a coordinate change of this form.
In either case, there is a constant $c_0$ for which  $|\partial_y^oS_i(x,y)| > c_0$ throughout all $D_{ikl}$. This means there is a nonvanishing $S_{0o}y^o$ term in the Taylor
expansion $\sum_{\alpha,\beta}S_{\alpha\beta}x^{\alpha}y^{\beta}$ of $S_i(x,y)$. So $\alpha_i$, so the minimum $\alpha + M_i\beta$ for nonvanishing $S_{\alpha\beta}$ must
be at most $oM_i$. Since $p_i(x,y) = \sum_{\alpha + M_i\beta = \alpha_i}  S_{\alpha\beta}x^{\alpha}y^{\beta}$, the degree of $p_i(1,y)$ is therefore at most
${\alpha_i \over M_i}\leq o$. Denote this  maximum power by $n_i$. Next, note that in view of $(5.7)$ we have
$$ \partial_{yy} (S_i(x,x^{M_i}(y + r_{ik}))) = x^{\alpha_i} \partial_{yy} p_i(1,y + r_{ik}) +  x^{\alpha_i + \xi} \partial_{yy} q(x,y + r_{ik}) $$
By the above discussion, $p_i(1,y + r_{ik})$ is a polynomial of degree $\leq n_i$. When one does the coordinate change $(x,y) 
\rightarrow (x, \pm y + \phi_{ikl}(x))$ transferring into the new coordinates, $ x^{\alpha_i} \partial_{yy} p_i(1,y + r_{ik})$ becomes 
$ x^{\alpha_i} \partial_{yy} p_i(1,\pm y + r_{ik})$ , while $x^{\alpha_i + \xi} \partial_{yy} q(x,y + r_{ik})$ transforms into some function of the form 
$x^{\alpha_i + \xi} s(x,y)$. Thus the sum of the terms of the Taylor expansion of $Q_{ikl}(x,y)$ with minimal $x$-power, given by 
$x^{\alpha_i} \partial_{yy} p_i(1,\pm y + r_{ik})$, can be written in the form $Cx^{\alpha_i}y^{n_{ikl}} + o(x^{\alpha_i}y^{n_{ikl}})$,
where $n_{ikl} \leq n_i -2 \leq o - 2$. Note that for the $D_{ikl}'$ for which $|y| > x^{\eta}$, since $\alpha_{ikl} = \alpha_i$ and $\beta_{ikl} = n_{ikl}$, the lemma follows from this. For the rest of the $D_{ikl}'$, we argue as follows.

Suppose a wedge $\{(x,y): 0 < x < e, \,\,c_1 x^m < y < c_2 x^m\}$ is contained in the domain $D_{ikl}'$, where $c_2 > c_1 > 0$. If one changes 
coordinates on this wedge, turning the former $(x,y)$ into $(x, x^my)$, then the wedge becomes the rectangle $K = \{(x,y): 0 < x < e,\,\, c_1 < y < c_2\}$, and the fact that $Q_{ikl}(x,y) \sim x^{\alpha_{ikl}}y^{\beta_{ikl}}$ on $D_{ikl}'$ implies that $Q_{ikl}(x,x^m y) \sim x^{\alpha_{ikl} + m \beta_{ikl}}y^{\beta_{ikl}}$ on $K$. Thus the terms of the Taylor series of $Q_{ikl}(x,x^m y)$ have 
$x$-degree at least $\alpha_{ikl} + m \beta_{ikl}$. Since the term $x^{\alpha_i}y^{n_{ikl}}$ becomes $x^{\alpha_i + mn_{ikl}}y^{n_{ikl}}$ and is one of the terms of the  Taylor series of $Q_{ikl}(x,x^m y)$, we must have that 
$\alpha_i + mn_{ikl} \geq \alpha_{ikl} + m \beta_{ikl}$. 

So if $D_{ikl}'$ has upper boundary $H_{ikl}x^{M_{ikl}} + ...$ and lower boundary $h_{ikl}x^{m_{ikl}}+ ...$, we must have $\alpha_i + mn_{ikl} \geq \alpha_{ikl} + m \beta_{ikl}$ for all $m_{ikl} \geq m \geq M_{ikl}$. So for all $(x,y)$ satisfying $x^{M_{ikl}} \geq y \geq x^{m_{ikl}}, 0 < x < 1$, one has $x^{\alpha_i}y^{n_{ikl}}
\leq x^{\alpha_{ikl}}y^{\beta_{ikl}}$. Thus there is a constant $C$ such that on the entire domain $D_{ikl}'$ one has 
 $Cx^{\alpha_i}y^{n_{ikl}} \leq x^{\alpha_{ikl}}y^{\beta_{ikl}}$. Similarly, if $D_{ikl}'$ has upper boundary $H_{ikl}x^{M_{ikl}} + ...$ and lower boundary the $x$-axis (corresponding to the case where $\beta_{ikl} = 0$), then $\alpha_i + mn_{ikl} \geq \alpha_{ikl} + m \beta_{ikl}$ for all $m \geq M_{ikl}$. So one analogously has $x^{\alpha_i}y^{mn_{ikl}}
\leq x^{\alpha_{ikl}}y^{\beta_{ikl}}$ for $y  < x^{M_{ikl}}$, $0 < x < 1$, and therefore  $Cx^{\alpha_i}y^{n_{ikl}} \leq x^{\alpha_{ikl}}y^{\beta_{ikl}}$ on all of $D_{ikl}'$ once again. Since $n_{ikl} \leq o-2$, we conclude that in either situation, one has
$Cx^{\alpha_i}y^{o-2} \leq x^{\alpha_{ikl}}y^{\beta_{ikl}}$ on  $D_{ikl}'$ and we are done with the proof of the lemma.

\noindent {\bf Estimates when $y$ is near a zero of $\partial_y p_i(1,y)$ of order greater than 1}.

We will bound the contribution to $T_{ij}^{N,2}(\lambda)$ coming from the integral $(5.9b)$ over the 
domain $D_{ikl}$ and add over all (finitely many) $k$ and $l$ to obtain the necessary estimates
for $|T_{ij}^{N,2}(\lambda)|$. Denote this integral over $D_{ikl}$ by $J_{ijkl}$. Performing the coordinate change $\phi_{ikl}$ we see that  $J_{ijkl}$ is given by
$$J_{ijkl}=  \int_{\{(x,y) \in  D_{ikl}':\, x \in [2^{-j-1}, 2^{-j}]\}}e^{i\lambda_1 x^{\alpha_i}S_{ikl}(x,y) + i\lambda_4x + i\lambda_3(x^{M_i}( y  + r_{ik} + \phi_{ikl}(x))+ \psi_i(x))}$$
$$\times  g(x^{\alpha_i}S_{ikl}(x,y))K_i (x,x^{M_i}(y + r_{ik} +  \phi_{ikl}(x)))\rho(2^jN(x - x_0))x^{M_i}\,dx\,dy \eqno (5.21)$$
(Without losing generality we are using $y + r_{ik}$ rather than $\pm y + r_{ik}$ to simplify the notation here.) We divide $(5.21)$ dyadically in the $y$ variable, writing $J_{ijkl} = \sum_m J_{ijklm}$, where 
$$J_{ijklm}=  \int_{([2^{-j-1}, 2^{-j}] \times [2^{-m-1}, 2^{-m}]) \cap D_{ikl}'}e^{i\lambda_1 x^{\alpha_i}S_{ikl}(x,y) + i\lambda_4x + i\lambda_3(x^{M_i} (y  + r_{ik} + \phi_{ikl}(x)) + \psi_i(x))} $$
$$\times g(x^{\alpha_i}S_{ikl}(x,y)) K_i (x,x^{M_i}(y + r_{ik} + \phi_{ikl}(x)))\rho(2^jN(x - x_0))x^{M_i}\,dx\,dy \eqno (5.22)$$
The estimates we need will be obtained by applying Lemma 3.1 to $(5.22)$ twice, once in the $x$ direction and once in the $y$ direction, and then taking the better of the two estimates thereby obtained. For the moment, we assume $M_i > s_i$ in Theorem 2.1 b) and will deal with the $M_i = s_i$ case afterwards.

We proceed to the application of Lemma 3.1 in the $x$ derivative, which will be used for second or third derivatives.  We examine the phase function in $(5.22)$. First, note that for some $\sigma > 0$ one has 
$$\partial_{xx}(x^{\alpha_i}S_{ikl}(x,y))  = \lambda_1\alpha_i(\alpha_i - 1)x^{\alpha_i - 2} S_{ikl}(0,y) + O(|\lambda_1|x^{\alpha_i - 2 + \sigma}) \eqno (5.23)$$
Next, note that since $\psi_i(x) = l_ix^{s_i} + O(x^{s_i + \sigma})$ with $l_i \neq 0$ for some $\sigma > 0$ (which we can take to be the same
as the previous $\sigma$), and since for the moment we are assuming that $M_i > s_i$, we analogously have 
$$\partial_{xx}(\lambda_3(x^{M_i} (y  + r_{ik} + \phi_{ikl}(x)) + \psi_i(x)))= \lambda_3l_i s_i(s_i - 1)x^{s_i - 2} + O(|\lambda_3| x^{s_i - 2 + \sigma}) \eqno (5.24)$$
Putting $(5.23)$ and $(5.24)$ together, if $P_{ikl}(x,y)$ denotes the phase function in $(5.22)$ we see that
$$\partial_{xx}P_{ikl}(x,y) =  \lambda_1\alpha_i(\alpha_i - 1)x^{\alpha_i - 2} S_{ikl}(0,y) + \lambda_3l_i s_i(s_i - 1)x^{s_i - 2}$$
$$+ 
O(|\lambda_1|x^{\alpha_i - 2 + \sigma} + |\lambda_3| x^{s_i - 2 + \sigma}) \eqno (5.25)$$
Analogously, one has 
$$\partial_{xxx}P_{ikl}(x,y) =  \lambda_1\alpha_i(\alpha_i - 1)(\alpha_i - 2)x^{\alpha_i - 3} S_{ikl}(0,y) + \lambda_3l_i s_i(s_i - 1)(s_i - 2) x^{s_i - 3}$$ 
$$+ O(|\lambda_1|x^{\alpha_i - 3 + \sigma} + |\lambda_3| x^{s_i - 3 + \sigma}) \eqno (5.26)$$
Since $\alpha_i > s_i > 1$,  the determinant of the matrix $M$ with rows $(\alpha_i(\alpha_i - 1), l_i s_i(s_i - 1))$ and $(\alpha_i(\alpha_i - 1)(\alpha_i - 2), l_i s_i(s_i - 1)(s_i - 2))$, given by $\alpha_i(\alpha_i - 1)l_i s_i(s_i - 1)(s_i - \alpha_i)$, is nonzero. As a result, there 
is a $c_0 > 0$ such that if $v$ denotes the vector $(\lambda_1x_0^{\alpha_i - 2} S_{ikl}(0,y), \lambda_3x_0^{s_i - 2})$
 where $x_0$ is as in the $\rho(2^jN(x - x_0))$ factor in $(5.22)$, then we have 
$$|M v| > c_0 |v| $$
Thus at least one of the components of $M v$ has magnitude at least ${c_0 \over 2} |v|$. If it is the first component, we have
$$|\lambda_1\alpha_i(\alpha_i - 1)x_0^{\alpha_i - 2} S_{ikl}(0,y) + \lambda_3l_i s_i(s_i - 1)x_0^{s_i - 2}| $$
$$ > {c_0 \over 2} (|\lambda_1\alpha_i(\alpha_i - 1)x_0^{\alpha_i - 2} S_{ikl}(0,y)| + |\lambda_3l_i s_i(s_i - 1)x_0^{s_i - 2}|) \eqno (5.27a)$$
If it is the second component, we have
$$|\lambda_1\alpha_i(\alpha_i - 1)(\alpha_i - 2) x_0^{\alpha_i - 3} S_{ikl}(0,y) + \lambda_3l_i s_i(s_i - 1)(s_i - 2)x_0^{s_i - 3}|$$
$$ > {c_0  \over 2} (|\lambda_1\alpha_i(\alpha_i - 1)(\alpha_i - 2)x_0^{\alpha_i - 3} S_{ikl}(0,y)| + |\lambda_3l_i s_i(s_i - 1)(s_i - 2)x_0^{s_i - 3}|) \eqno (5.27b)$$
Assuming the parameter $N$ in $(5.22)$ was chosen sufficiently large, for a given $y$ equation $(5.27a)$ or $(5.27b)$ will not just hold at $x = x_0$, but for all $x$  in the domain of the integrand of $(5.22)$. Furthermore, if $x$ is sufficiently small,
the error terms in $(5.25)$ or $(5.26)$ respectively will be of magnitude at most half that of the right hand side of $(5.27a)$ or
$(5.27b)$ respectively. (Here we implicitly use that $|S_{ikl}(0,y)|$ is bounded below over $y \in [0,H_i]$, but recall $S_{ikl}(0,y)$ never vanishes 
on $[0,H_i]$ since $S_i(x,y) \sim x^{\alpha_i}$). Also, one technical point worth mentioning here: If $s_i = 2$, the second term
of the error term of $(5.26)$ will not be small in comparison to $|\lambda_3l_i s_i(s_i - 1)(s_i - 2)x_0^{s_i - 3}|)$, but it will be small
in comparison with $|\lambda_1\alpha_i(\alpha_i - 1)(\alpha_i - 2)x_0^{\alpha_i - 3} S_{ikl}(0,y)| $ if $(5.27a)$ does not hold since
$\alpha_i > s_i = 2$. Analogous considerations would hold if $\alpha_i = 2$ since then $s_i \neq 2$, but one can show that $\alpha_i =2$ never occurs.

Thus $(5.27a)$ or $(5.27b)$ imply that for a fixed $y$, for all $x$ in the 
domain of integration of $(5.22)$ we have at least one of the following equations holding.
$$|\partial_{xx}P_{ikl}(x,y)| > {c_0 \over 4}(|\lambda_1\alpha_i(\alpha_i - 1)x^{\alpha_i - 2} S_{ikl}(0,y)| + |\lambda_3l_i s_i(s_i - 1)x^{s_i - 2}|) \eqno (5.28a)$$
$$|\partial_{xxx}P_{ikl}(x,y)|> {c_0 \over 4} (|\lambda_1\alpha_i(\alpha_i - 1)(\alpha_i - 2)x^{\alpha_i - 3} S_{ikl}(0,y)| + |\lambda_3l_i s_i(s_i - 1)(s_i - 2)x^{s_i - 3}|) \eqno (5.28b)$$
What is relevant for us is that for some constant $c_1 > 0$, for each fixed $y$ we have 
$$ |\partial_{xx}P_{ikl}(x,y)|  > c_1 |\lambda_1|x^{\alpha_i - 2} \,\,\,\,\,\,\,\,\,\,\,\,\,\,\,\,\,\,\,\,{\rm or} \,\,\,\,\,\,\,\,\,\,\,\,\,\,\,\,\,\,\,\,
|\partial_{xxx}P_{ikl}(x,y)|>  c_1 |\lambda_1|x^{\alpha_i - 3} \eqno (5.29a)$$
Equivalently, for some constant $c_2$ we have
$$ |\partial_{xx}P_{ikl}(x,y)|  > c_2 |\lambda_1|2^{-j(\alpha_i - 2)} \,\,\,\,\,\,\,\,\,\,\,\,\,\,\,\,\,\,\,\,{\rm or} \,\,\,\,\,\,\,\,\,\,\,\,\,\,\,\,\,\,\,\,
|\partial_{xxx}P_{ikl}(x,y)|>  c_2 |\lambda_1|2^{-j(\alpha_i - 3)} \eqno (5.29b)$$
Equation $(5.29b)$ is what we will use to apply Lemma 3.1 on the $x$-integral of $(5.22)$. 

Note that in the above we assumed 
$M_i > s_i$, and next we will show that $(5.29b)$ also holds when $M_i = s_i$. By Theorem 2.1, $M_i = s_i$ can only occur if $\psi_i(x) 
= l_i x^{s_i}$ for some $l_i \neq 0$. So the term $\lambda_3(x^{M_i} (y  + r_{ik} + \phi_{ikl}(x)) + \psi_i(x))$ in the exponent of $(5.22)$
can be written as $\lambda_3 x^{M_i}(y + r_{ik} +  l_i + \phi_{ikl}(x))$. Recall by $(4.10)$ that $\lambda_3 = -{d_i\alpha_i(\alpha_i - 1) \over l_is_i(s_i - 1)}x_0^{\alpha_i - s_i} = -{d_i\alpha_i(\alpha_i - 1) \over l_is_i(s_i - 1)}x_0^{\alpha_i - M_i}$. So there is some constant $c$ such that the term is equal to
$\lambda_1 cx^{M_i}x_0^{\alpha_i - M_i}((y + r_{ik} +  l_i) + \phi_{ikl}(x))$. This means that if $|y + r_{ik} +  l_i|$ is sufficiently
small, the second $x$ derivative of this term of the phase in $(5.22)$ is negligible in comparison to $\lambda_1\alpha_i(\alpha_i - 1)x^{\alpha_i - 2} S_{ikl}(0,y)$ and therefore the left-hand side of $(5.29b)$ holds. On the other hand if $|y + r_{ik} +  l_i|$ is
not small, the argument used to show $(5.29b)$ applies once again, so once again $(5.29b)$ holds.

We now apply Lemma 3.1  on the $x$-integral of $(5.22)$ using $(5.29b)$. To do this we also have to bound the $A(x,y) =  g(x^{\alpha_i}S_{ikl}(x,y)) K_i (x,x^{M_i}(y + r_{ik} + \phi_{ikl}(x)))\rho(2^jN(x - x_0))x^{M_i}$ and the integral of its 
$x$ derivative for a given $y$. Using $(1.2a)-(1.2b)$ analogously to $(3.32)$ we have
$$|A(x,y)| \leq (2^{-j\alpha_i\alpha})(2^{-j\beta})(2^{-jM_i}) \eqno (5.30)$$
I claim that by taking an $x$ derivative of $A(x,y)$ one gains an additional factor bounded by $C{2^j}$. To show this,
 by the product rule for derivatives it suffices to show that by differentiating each of the factors of $A(x,y)$ one gains at most $C{1 \over x}$.
For the $x^{M_i}$ factor this is obvious. The same is true for the $\rho(2^jN(x - x_0))$ factor since $x \sim 2^{-j}$ and $N$ 
is a constant. For the $g(x^{\alpha_i}S_{ikl}(x,y))$ factor we have
$$|\partial_x g(x^{\alpha_i}S_{ikl}(x,y))| \leq Cx^{\alpha_i - 1} |g'(x^{\alpha_i}S_{ikl}(x,y))| \eqno (5.31)$$
By $(1.2a)$ this is bounded by 
$$C x^{\alpha_i - 1} (x^{\alpha_i})^{\alpha - 1}  = Cx^{\alpha_i \alpha - 1} \eqno (5.32)$$
Lastly, for $  K_i (x,x^{M_i}(y + r_{ik} + \phi_{ikl}(x)))| $ what we need follows from $(1.2b)$ and the fact that $M_i \geq 1$. So we conclude
that we have
$$|\partial_x A(x,y)| \leq (2^{-j\alpha_i\alpha})(2^{-j\beta})(2^{-jM_i}) (2^j)\eqno (5.33)$$
This is the estimate we will use for $|\partial_x A(x,y)| $ in applying Lemma 3.1.

We now apply Lemma 3.1 for fixed $y$ in $(5.22)$, using $(5.29b)$ on the phase and $(5.32)$, and $(5.33)$ on $A(x,y)$. Afterwards, we integrate the result in $y$. We obtain
$$|J_{ijklm}| \leq C (2^{-j-m}) (2^{-j\alpha_i\alpha})(2^{-j\beta})(2^{-jM_i})\max \bigg({1 \over (|\lambda_1|2^{-j\alpha_i})^{1 \over 2}}, {1 \over (|\lambda_1|2^{-j\alpha_i})^{1 \over 3}}\bigg) \eqno (5.34)$$
Simply by taking absolute values and integrating in $(5.22)$, one has
$$|J_{ijklm}| \leq C (2^{-m}) (2^{-j\alpha_i\alpha})(2^{-j\beta})(2^{-jM_i})(2^{-j}) \eqno (5.35)$$
Note that the left hand side of the maximum in $(5.34)$ is greater than the right hand side if and only if $|\lambda_1|2^{-j\alpha_i} < 1$, in which case $(5.35)$ gives a better estimate anyhow. Thus $(5.34)$ and $(5.35)$ combine into
$$|J_{ijklm}| \leq C (2^{-j-m}) (2^{-j\alpha_i\alpha})(2^{-j\beta})(2^{-jM_i})\min\bigg(1, {1 \over (|\lambda_1|2^{-j\alpha_i})^{1 \over 3}}\bigg)\eqno (5.36)$$
We now examine the estimates obtained by applying Lemma 3.1 to $(5.22)$ in the $y$-direction. Note that
$$\partial_{yy} P_{ikl}(x,y) = \lambda_1 \partial_{yy}(x^{\alpha_i}S_{ikl}(x,y)) \eqno (5.37)$$
$$=  \lambda_1{Q}_{ikl}(x,y) \eqno (5.38)$$
Recall that $Q_{ikl}(x,y)$ is comparable to the monomial $x^{\alpha_{ikl}}y^{\beta_{ikl}}$ on $D_{ikl}$ with $\alpha_{ikl} \geq \alpha_i$. So $(5.38)$ implies that
$$|\partial_{yy} P_{ikl}(x,y)| > C|\lambda_1|x^{\alpha_{ikl}}y^{\beta_{ikl}} \eqno (5.39)$$
So on the support of the integrand of $(5.22)$ one has
$$|\partial_{yy} P_{ikl}(x,y)| > C|\lambda_1|2^{-j\alpha_{ikl} - m\beta_{ikl}} \eqno (5.40)$$
In $(5.22)$, we now apply Lemma 3.1 in the $y$ direction, using $(5.40)$, $(5.12)$, and $(5.15)$, and integrate the result in $x$.
(We can still use $(5.12)$ and $(5.15)$ here due to the form of the coordinate change $(x,y) \rightarrow (x,y + \phi_{ikl}(x))$ ).
We get that
$$|J_{ijklm}|  \leq C (2^{-j\alpha_i\alpha})( 2^{-j\beta}) (2^{-j M_i}){1 \over (|\lambda_1|2^{-j\alpha_{ikl} -m\beta_{ikl}})^{1 \over 2}}(2^{-j})  \eqno (5.41)$$
Equations $(5.36)$ and $(5.41)$ can be combined into a single estimate:
$$|J_{ijklm}| \leq C (2^{-j-m}) (2^{-j\alpha_i\alpha -j\beta})(2^{-jM_i})\min\bigg(1, {1 \over (|\lambda_1|2^{-j\alpha_i})^{1 \over 3}}, {1 \over (|\lambda_1|2^{-j\alpha_{ikl} -m(\beta_{ikl} + 2)})^{1 \over 2}}\bigg)\eqno (5.42)$$
This is equivalent to 
$$|J_{ijklm}| \leq C\int_{[2^{-j-1},2^{-j}] \times [2^{-m-1},2^{-m}]} x^{\alpha_i\alpha + \beta + M_i}\min\bigg(1, {1 \over (|\lambda_1|x^{\alpha_i})^{1 \over 3}}, {1 \over (|\lambda_1| x^{\alpha_{ikl}}y^{\beta_{ikl}+ 2})^{1 \over 2}}\bigg)\eqno (5.43)$$
In view of the shapes of the $D_{ikl}'$, adding this over all $m$ therefore gives the following.
$$|J_{ijkl}| \leq C \int_{\{(x,y) \in  D_{ikl}': x \in [2^{-j-1}, 2^{-j}]\}}x^{\alpha_i\alpha + \beta +  M_i}\min\bigg(1, {1 \over (|\lambda_1|x^{\alpha_i})^{1 \over 3}}, {1 \over (|\lambda_1| x^{\alpha_{ikl}}y^{\beta_{ikl}+ 2})^{1 \over 2}}\bigg)\eqno (5.44)$$
By Lemma 5.1, we therefore have
$$|J_{ijkl}| \leq C \int_{\{(x,y) \in  D_{ikl}': x \in [2^{-j-1}, 2^{-j}]\}}x^{\alpha_i\alpha + \beta + M_i}\min\bigg(1, {1 \over (|\lambda_1|x^{\alpha_i})^{1 \over 3}}, {1 \over (|\lambda_1|x^{\alpha_i}y^o)^{1 \over 2}}\bigg)\,dx\,dy\eqno (5.45)$$
We fix $x$ and focus on the $y$ integral of $(5.45)$, which is at most
$$x^{\alpha_i\alpha + \beta + M_i}\int_0^1 \min\bigg({1 \over (|\lambda_1|x^{\alpha_i})^{1 \over 3}}, {1 \over (|\lambda_1|x^{\alpha_i}y^o)^{1 \over 2}}\bigg)\,dx\,dy\eqno (5.46)$$
The quantities $(|\lambda_1|x^{\alpha_i})^{1 \over 3}$ and $(|\lambda_1|x^{\alpha_i}y^o)^{1 \over 2}$ are equal when
$(|\lambda_1|x^{\alpha_i})^2 = (|\lambda_1|x^{\alpha_i}y^o)^3$, in other words when $y^{3o} = |\lambda_1|^{-1}x^{-\alpha_i}$ or $y =   |\lambda_1|^{-{1 \over 3o}}x^{-{\alpha_i \over 3o}}$. The integrand in $(5.46)$ is constant for 
$y$ smaller than this value, and decreases like $y^{-{o \over 2}}$ for $y$ larger than this value. Recalling that $o \geq 3$ in the situation at hand since $\partial_y p_i(1,y)$ has a zero of order at least two at $r_{ik}$, $(5.46)$ is therefore bounded
by a constant times $x^{\alpha_i\alpha + \beta + M_i} \times|\lambda_1|^{-{1 \over 3o}}x^{-{\alpha_i \over 3o}} \times {1 \over (|\lambda_1|x^{\alpha_i})^{1 \over 3}}
= {x^{\alpha_i\alpha + \beta + M_i} \over (|\lambda_1|x^{\alpha_i})^{{1 \over 3} + {1 \over 3o}}}$. Integrating this in $x$ we
therefore have
$$|J_{ijkl}| \leq C2^{-j\alpha_i\alpha - j\beta - jM_i - j} {1 \over (|\lambda_1|2^{-j\alpha_i})^{{1 \over 3}+ {1 \over 3o}} }
\eqno (5.47)$$
By just taking the $1$ in the minimum of $(5.45)$ and integrating one gets
$$|J_{ijkl}| \leq C2^{-j\alpha_i\alpha - j\beta - jM_i - j} \eqno (5.48)$$
So combining $(5.47)$ and $(5.48)$ we see that
$$|J_{ijkl}| \leq C2^{-j\alpha_i\alpha - j\beta - jM_i - j}\min\bigg(1, {1 \over (|\lambda_1|2^{-j\alpha_i})^{{1 \over 3}+ {1 \over 3o}} }\bigg)\eqno (5.49)$$
Given that $D_i'$ was of $y$-width $\sim x^{M_i}$ for a given $x$, analogously to $(5.18b)$ and $(3.23)$  equation $(5.49)$ implies
$$|J_{ijkl}| \leq C \int_{\{(x,y) \in D_i':\,x \in [2^{-j-1},2^{-j}]\}}\min\bigg(1, {1 \over |\lambda_1|^{{1 \over 3}+ {1 \over 3o}}|S_i(x,y)|^{{1 \over 3}+ {1 \over 3o}}}\bigg)\,d\mu_{\alpha,\beta} \eqno (5.50)$$

Note that the right-hand side of $(5.50)$ is independent of $k$ and $l$. So if one adds over all $k$ and $l$, the result is bounded by
the right-hand side of $(5.50)$.
But the sum over all $k$ and $l$ of $J_{ijkl}$ is exactly $J_{ij}$, the contribution to $T_{ij}^{N,2}(\lambda)$ in the integral $(5.9b)$ coming from the sets $[2^{-j-1},2^{-j}] \times ([0,H_i] \cap I_k)$ where $p'(1,y)$ has a zero of order at least 2 at $y = r_{ik}$. Note that since ${1 \over 3} + {1 \over 3o} \leq {1 \over 2}$, the integrand in $(5.50)$ is at least as large as that of $(5.18b)$, which by $(5.20)$ is the bound we have for the contribution to $T_{ij}^{N,2}(\lambda)$ in the integral $(5.9b)$ coming from the sets $[2^{-j-1},2^{-j}] \times ([0,H_i] \cap I_k)$ where $p'(1,y)$ has a zero of order $1$ at $y = r_{ik}$. So adding over all $k$, we have the succinct statement that 
$$|T_{ij}^{N,2}(\lambda)| \leq C \int_{\{(x,y) \in D_i':\,x \in [2^{-j-1},2^{-j}]\}}\min\bigg(1, {1 \over |\lambda_1|^{{1 \over 3}+ {1 \over 3o}}|S_i(x,y)|^{{1 \over 3}+ {1 \over 3o}}}\bigg)\,d\mu_{\alpha,\beta} \eqno (5.51)$$
Thus to complete the proof of Theorem 1.1, it suffices to show that the right-hand side of $(5.51)$ is bounded by the right-hand
sides of $(1.8a)-(1.8c)$. But $(5.51)$ is the same as $(3.24)$, and the steps from $(3.24)$ to $(3.31)$ give $(1.8a)-(1.8c)$ 
exactly as before. This completes the proof of Theorem 1.1.

\noindent {\bf 6. Sharpness of Theorem 1.1a) when $\beta = 0$.}

To see why when $\beta = 0$ the uniform estimates given by Theorem 1.1a) are sharp and that $(\delta,d) = (\alpha + \delta_0,d_0)$, where $(\delta_0,d_0)$ are the $(\delta,d)$ of the smooth case (with $\alpha = \beta = 0$), we use some facts concerning the asymptotics of sublevel set measures and their connection to oscillatory integrals that follow from two-dimensional resolution of singularities. We refer to [AGV] ch 7 for more information.
If $E_r$ denotes the disk $\{(x,y): x^2 + y^2 < r^2\}$, then if $r$ is sufficiently small by resolution of singularities one has an asymptotic expansion
$$|{\{(x,y) \in E_r :  |S(x,y)| < t\}}| =  D_r t^{\delta_0} |\ln t|^{d_0} + o (t^{\delta_0} |\ln t|^{d_0})\eqno (6.1)$$
Here $D_r \neq 0$. The terms of the $o(t^{\delta_0} |\ln t|^{d_0})$ part of the asymptotics are of the form $ct^a |\ln t|^b$, where 
$b = 0$ or $1$ and $a$ is a rational number. The set of all such possible $a$ are a subset of an arithmetic progression whose smallest
value is greater than $\delta_0$. 

\noindent Next, one has
$$ \int_{\{(x,y) \in E_r :  |S(x,y)| < \epsilon\}}|S(x,y)|^{ \alpha}\,dx\,dy = \int_0^{\epsilon} t^{\alpha} {\partial \over \partial t} 
\bigg(D_r t^{\delta_0} |\ln t|^{d_0} + o (t^{\delta_0} |\ln t|^{d_0})\bigg)\,dt \eqno (6.2)$$
We can assume that $\alpha  > -\delta_0$ since $|S(x,y)|^{\alpha}$ is not locally integrable on a  neighborhood of the origin if 
$\alpha \leq -\delta_0$. So $(6.2)$ is of the form $C\epsilon^{\alpha + \delta_0} |\ln \epsilon|^{d_0}$ plus a smaller error term. So the $(\delta,d)$ in $(1.6)$ is given by $(\alpha + \delta_0, d_0)$.  Note that in particular, by the form  of $(1.8a)$, one has
$\delta < {1 \over 3} + {1 \over o} < 1$. So if Theorem 1.1a) holds then $\alpha + \delta_0 < 1$.

If $\phi(x,y)$ is a smooth function
supported in a sufficiently small neighborhood of the origin then we have analogous asymptotics of the following form as $t \rightarrow 0^+$.
$$\int_{\{(x,y): 0 < S(x,y) < t\}} \phi(x,y)\,dx\,dy = C_{\phi} t^{\delta_0} |\ln t|^{d_0} + o (t^{\delta_0} |\ln t|^{d_0})\eqno (6.3a)$$
$$\int_{\{(x,y): 0 > S(x,y) > -t\}} \phi(x,y)\,dx\,dy = C_{\phi}' t^{\delta_0} |\ln t|^{d_0} + o (t^{\delta_0} |\ln t|^{d_0})\eqno (6.3b)$$
At least one of $C_{\phi}$ and $C_{\phi}'$ will be 
nonzero if $\phi(x,y)$ is nonnegative with $\phi(0,0) > 0$.  We further have that
$$\int|S(x,y)|^{\alpha} e^{i\lambda_1 S(x,y)}\phi(x,y)\,dx\,dy = \int_0^1 t^{\alpha} e^{i\lambda_1 t} {\partial \over \partial t} \bigg(\int_{\{(x,y): 0 < S(x,y) < t\}} \phi(x,y)\,dx\,dy\bigg)\,dt$$
$$+ \int_0^1 t^{\alpha} e^{-i\lambda_1 t} {\partial \over \partial t} \bigg(\int_{\{(x,y): -t < S(x,y) < 0\}} \phi(x,y)\,dx\,dy\bigg)\,dt \eqno (6.4)$$
$$= \int_0^1 t^{\alpha} e^{i\lambda_1 t}( C_{\phi} \delta_0 t^{\delta_0 - 1} |\ln t|^{d_0} + o (t^{\delta_0 - 1} |\ln t|^{d_0}))\,dt$$
$$+ \int_0^1 t^{\alpha} e^{-i\lambda_1 t}( C_{\phi}' \delta_0 t^{\delta_0 - 1} |\ln t|^{d_0} + o (t^{\delta_0 - 1} |\ln t|^{d_0}))\,dt
\eqno (6.5)$$
$$ =  C_{\phi} \delta_0\int_0^{1}  e^{i\lambda_1 t}t^{\delta_0 + \alpha - 1} |\ln t|^{d_0}\,dt  + C_{\phi}' \delta_0\int_0^{1}  e^{-i\lambda_1 t}t^{\delta_0 + \alpha - 1} |\ln t|^{d_0}\,dt $$
$$+ \int_0^{1} o (t^{\delta_0 + \alpha - 1} |\ln t|^{d_0})\,dt
\eqno (6.6)$$
As described above, we have that $0 < \alpha + \delta_0 < 1$ whenever Theorem 1.1a) holds with $\beta = 0$. In this case, the first two terms in $(6.6)$ can be computed directly to be of the form $C\lambda_1^{-\alpha - \delta_0}(\ln \lambda_1)^{d_0}$ ($C \neq 0$) plus a faster-decaying term as $\lambda_1 \rightarrow \infty$. These two main
terms will not cancel out because that calculation also reveals these main terms will not be real multiples of each other.

One can expand the $o (t^{\delta_0 + \alpha - 1} |\ln t|^{d_0})$ term in $(6.4)$ to any finite number of terms, and
one can do a similar integration on each term obtained to get terms decaying faster than $C\lambda_1^{-\alpha - \delta_0}(\ln \lambda_1)^{d_0}$.  The integral corresponding to the 
error term in such an expansion can be bounded using integration by parts, and the decay rate of the error term increases 
indefinitely with the number of terms
in the expansion. Thus for some $C \neq 0$, $\int|S(x,y)|^{\alpha} e^{i\lambda_1 S(x,y)}\phi(x,y)\,dx\,dy$ is equal to $C\lambda_1^{-\alpha - \delta_0}(\ln \lambda_1)^{d_0}$ plus a term that decays faster as $\lambda_1 \rightarrow \infty$. Since $(\delta, d) = (\delta_0 + \alpha, d_0)$, the estimate given by Theorem 1.1a) is therefore seen to sharp by letting $g(z) = |z|^{\alpha}$ and $K(x,y) = \phi(x,y)$ be a nonnegative function with $\phi(0,0) > 0$. 

\noindent We conclude that the estimates given by Theorem 1.1a) are sharp whenever $\beta = 0$.

\noindent {\bf 7. References.}

\noindent [AGV] V. Arnold, S. Gusein-Zade, A. Varchenko, {\it Singularities of differentiable maps},
Volume II, Birkhauser, Basel, 1988. \parskip = 4pt\baselineskip = 4pt

\noindent [CKaN] K. Cho, J. Kamimoto, T. Nose, {\it Asymptotic analysis of oscillatory integrals via the Newton polyhedra of the phase and the amplitude},  J. Math. Soc. Japan, {\bf 65} (2013), no. 2, 521-562.

\noindent [CoMa] M. Cowling, G. Mauceri, {\it Oscillatory integrals and Fourier transforms of surface-carried measures}, Trans. Amer.
Math. Soc. {\bf 304} (1987), no. 1, 53-68.

\noindent [D] J.J. Duistermaat, {\it Oscillatory integrals, Lagrange immersions, and unfolding of singularities}, Comm. Pure Appl.
Math., {\bf 27} (1974), 207-281.

\noindent [G1] M. Greenblatt, {\it Estimates for Fourier transforms of surface measures in $\R^3$ with PDE applications},
submitted.

\noindent [G2] M. Greenblatt, {\it Resolution of singularities in two dimensions and the stability of integrals}, Adv. Math., 
{\bf 226} no. 2 (2011), 1772-1802.

\noindent [G3] M. Greenblatt, {\it The asymptotic behavior of degenerate oscillatory integrals in two
dimensions}, J. Funct. Anal. {\bf 257} (2009), no. 6, 1759-1798.

\noindent [IKeM1] I. Ikromov, M. Kempe, and D. M\"uller, {\it Damped oscillatory integrals and boundedness of
maximal operators associated to mixed homogeneous hypersurfaces} (English summary) Duke Math. J. {\bf 126} 
(2005), no. 3, 471--490.

\noindent [IKeM2] I. Ikromov, M. Kempe, and D. M\"uller, {\it Estimates for maximal functions associated
to hypersurfaces in $R^3$ and related problems of harmonic analysis}, Acta Math. {\bf 204} (2010), no. 2,
151--271.

\noindent [IM1] I. Ikromov, D. M\"uller, {\it On adapted coordinate systems}, Trans. AMS, {\bf 363} (2011), 2821-2848.

\noindent [IM2] I. Ikromov, D. M\"uller, {\it  Uniform estimates for the Fourier transform of surface-carried measures in
$\R^3$ and an application to Fourier restriction}, J. Fourier Anal. Appl, {\bf 17} (2011), no. 6, 1292-1332.

\noindent [IoSa1] A. Iosevich, E. Sawyer, {\it Oscillatory integrals and maximal averages over homogeneous
surfaces}, Duke Math. J. {\bf 82} no. 1 (1996), 103-141.

\noindent [IoSa2] A. Iosevich, E. Sawyer, {\it Maximal averages over surfaces},  Adv. Math. {\bf 132} 
(1997), no. 1, 46--119.

\noindent [KaN] J. Kamimoto, T. Nose, {\it Toric resolution of singularities in a certain class of $C^{\infty}$ functions and asymptotic analysis of oscillatory integrals}, preprint.

\noindent [K1] V. N. Karpushkin, {\it A theorem concerning uniform estimates of oscillatory integrals when
the phase is a function of two variables}, J. Soviet Math. {\bf 35} (1986), 2809-2826.

\noindent [K2] V. N. Karpushkin, {\it Uniform estimates of oscillatory integrals with parabolic or 
hyperbolic phases}, J. Soviet Math. {\bf 33} (1986), 1159-1188.

\noindent [Lic] B. Lichtin, {\it Uniform bounds for two variable real oscillatory integrals and singularities of mappings}, J. Reine Angew. Math. {\bf 611} (2007), 1–73.

\noindent [Lit] W. Littman, {\it Fourier transforms of surface-carried measures and differentiability of surface averages}, Bull. Amer.
Math. Soc., {\bf 69} (1963), 766-770.

\noindent [PS1] D. H. Phong, E. M. Stein, {\it The Newton polyhedron and
oscillatory integral operators}, Acta Mathematica {\bf 179} (1997), 107-152.

\noindent [PS2] D. H. Phong, E. M. Stein, {\it Damped oscillatory integral operators with analytic phases},  Adv. Math. {\bf 134} (1998), no. 1, 146-177. 

\noindent [PrY] M. Pramanik, C.W. Yang, {\it Decay estimates for weighted oscillatory integrals in $\R^2$},
Indiana Univ. Math. J., {\bf 53}  (2004), 613-645.

\noindent [SoS] C. Sogge and E. Stein, {\it Averages of functions over hypersurfaces in $R^n$}, Invent.
Math. {\bf 82} (1985), no. 3, 543--556.

\noindent [S] E. Stein, {\it Harmonic analysis; real-variable methods, orthogonality, and oscillatory 
integrals}, Princeton Mathematics Series Vol. 43, Princeton University Press, Princeton, NJ, 1993.

\noindent [V] A. N. Varchenko, {\it Newton polyhedra and estimates of oscillatory integrals}, Functional 
Anal. Appl. {\bf 18} (1976), no. 3, 175-196.

\line{}
\line{}

\noindent Department of Mathematics, Statistics, and Computer Science \hfill \break
\noindent University of Illinois at Chicago \hfill \break
\noindent 322 Science and Engineering Offices \hfill \break
\noindent 851 S. Morgan Street \hfill \break
\noindent Chicago, IL 60607-7045 \hfill \break
\noindent greenbla@uic.edu \hfill\break

\end